\documentclass[preprint,amsmath,amssymb,aps,pre,superscriptaddress,nofootinbib]{revtex4-1}

\usepackage{xcolor}
\usepackage{natbib}
\usepackage{graphicx}
\usepackage{dcolumn}
\usepackage{mathtools}
\usepackage{tikz}
\usepackage{appendix}
\usepackage{amsfonts}
\usepackage{amsmath}
\usepackage{amssymb}

\usetikzlibrary{calc,patterns,angles,quotes}

\begin{document}

\title{On the trajectory of the nonlinear pendulum:\\  Exact analytic solutions via power series.}

\author{W. Cade Reinberger}
\affiliation{School of Mathematical Sciences, Rochester Institute of Technology, Rochester, NY 14623, USA}
\author{Morgan S. Holland}
\affiliation{School of Mathematical Sciences, Rochester Institute of Technology, Rochester, NY 14623, USA}
\author{Nathaniel S. Barlow}
\affiliation{School of Mathematical Sciences, Rochester Institute of Technology, Rochester, NY 14623, USA}
\author{Steven J. Weinstein}
\affiliation{Department of Chemical Engineering, Rochester Institute of Technology, Rochester, NY 14623, USA}
\affiliation{School of Mathematical Sciences, Rochester Institute of Technology, Rochester, NY 14623, USA}

\date{\today}     

\begin{abstract}

We provide an exact infinite power series solution that describes the trajectory of a nonlinear simple pendulum undergoing librating and rotating motion for all time. Although the series coefficients were previously given in [V. Fair{\'e}n, V. L{\'o}pez, and L. Conde, \textit{Am. J. Phys} \textbf{56} (1), (1988), pp. 57?61], the series itself \--- as well as the optimal location about which an expansion should be chosen to assure series convergence and maximize the domain of convergence \--- was not examined, and is provided here.    By virtue of its representation as an elliptic function, the pendulum function has singularities that lie off of the real axis in the complex time plane.  This, in turn, imposes a radius of convergence on the physical problem in real time. By choosing the expansion point at the top of the trajectory, the power series converges all the way to the bottom of the trajectory without being affected by these singularities.  In constructing the series solution, we re-derive the coefficients using an alternative approach that generalizes to other nonlinear problems of mathematical physics.  Additionally, we provide an exact resummation of the pendulum series \--- Motivated by the asymptotic approximant method given in [Barlow et al., \textit{Q. J. Mech. Appl. Math.}, \textbf{70} (1) (2017), pp. 21-48] \--- that accelerates the series' convergence uniformly from the top to the bottom of the trajectory. We also provide an accelerated exact resummation of the infinite series representation for  the elliptic integral used in calculating the period of a pendulum's trajectory.   This allows one to preserve analyticity in the use of the period to extend the pendulum series for all time via symmetry.
\end{abstract}

\maketitle

\section{Introduction \label{sec:intro}} 

The pendulum equation is often presented in undergraduate science and math courses as a physically-relevant nonlinear ordinary differential equation (ODE) whose general solution cannot be written in terms of elementary functions.  The ODE is used to motivate phase-plane and stability analyses~\cite{Simmons,Strogatz} and, when linearized for small deflections, admits a classical harmonic oscillator solution.  It has been known for over a century that the analytic solution to the pendulum equation, expressed as angle as a function of time, may be written in terms of Jacobi elliptic functions~\cite{Helmholtz, Appell, Whittaker}.  The inverse solution (time as a function of angle) is expressed in terms of elliptic integrals and has been well documented for even longer; its application towards the evaluation of the period of a pendulum's trajectory was established alongside those of elliptic integrals in the 18$^\mathrm{th}$ century~\cite{Stillwell,Gray}. The singularities of the Jacobi elliptic functions (and thus the pendulum equation itself) were studied during the early development of the field of complex analysis~\cite{Gray} and are now well-documented~\cite{Byrd}.  The efficient representation of Jacobi elliptic functions via power series is limited by the radius of convergence imparted by such singularities.  While the use of power series solutions for linear ODEs is commonplace, the analogous treatment for their nonlinear counterparts is only found in scattered pockets of the literature.  In general, power series solutions to nonlinear ODEs have been thought by many (in the words of \cite{Agnew} p. 263) to `` play minor or nonexistent roles in the lives of most people''. It is perhaps not surprising that a power series solution \--- for the pendulum angle as a function of time \--- is virtually absent from the literature, especially since numerical algorithms for elliptic functions are readily available in many scientific computing packages.  An early attempt at finding the power series solution to the pendulum equation is given in~\cite{Agnew} (p. 264), where substitution of the initial conditions directly into the ODE and then successive (and objectively arduous) applications the chain rule leads one after four terms to (again in the author's words) ``feed the wastebasket and go to work on the next chapter''.     A few decades later, a methodology was developed~\cite{Fairen} to obtain the recursion for all coefficients of the pendulum series through a clever decomposition of a system of recursions. That said, the series of~\cite{Fairen} was implied as being approximate (by the paper's title), and a study of the radius of convergence of the series was not provided.  In fact, the solution itself is incomplete for that reason, as the initial conditions are left arbitrary, and convergence is not guaranteed. The aim of this work is to show that the power series representation of the pendulum function given by~\cite{Fairen} is, in fact, exact and analytic over the entire physical domain if one judiciously chooses the initial conditions (or equivalently, expansion point) along a given trajectory to assure convergence. Once a convergent series is constructed through this choice of initialization, the series may be re-summed exactly to accelerate convergence, and thus can be used as an alternative to elliptic functions.  It is relevant to note here that, although expansions for Jacobi elliptic functions are well-documented~\cite{Byrd}, the number of arithmetic steps needed needed to convert such expansions to a single power series for the pendulum solution is larger than those to directly find the power series solution of the original nonlinear ODE; we demonstrate this through the analysis in the current work. More generally, insights provided in this paper regarding power series solutions to nonlinear ODEs are relevant to problems that cannot be expressed by special functions such as Jacobi elliptic functions. 

%  The denominator of a Pad\'e serves to approximately incorporate the known (or unknown) singularities responsible for the divergence of the original series.  On the other hand, a Pad\'e has the undesired feature of introducing new singularities while also constrain.  While there are well-studied virtues of the Pad\'e and its ability to accelerate the convergence of a series~\cite{bender,BakerBook}, drawbacks include spurious singularities that it can inflict upon a problem and its rigid asymptotic behavior (away from the original expansion point) by virtue of the Pad\'e being a rational function~\cite{BakerGammel}.

In general, power series solutions to ODEs are limited by their radius of convergence as mentioned above.  In linear ODEs, the radius of convergence may be anticipated based on the distance between the expansion point and singularities in ODE coefficients.  In nonlinear ODE's, singularities may be spontaneous and cannot be anticipated by inspection.  As such, the radius of convergence imparted to power series solutions of nonlinear ODEs can restrict an analytic representation over the whole domain of interest \-- or, can lead to slow convergence properties.  Thus, power series solutions are often deficient, and with the prevalence of numerical techniques nowadays, often ignored as a viable solution technique.  However, techniques have been developed to accelerate convergence and even sum divergent series beyond the convergence radius (i.e., analytic continuation).  Approximant techniques, such as the well-known Pad\'e, have been utilized to this end~\cite{Boyd1997}.    A relatively new approach, \textit{asymptotic approximants}~\cite{Barlow:2017}, has been successfully applied to power series expansions to create highly accurate approximate solutions of nonlinear ODEs arising in many areas of mathematical physics~\cite{FS2020,SIR2020,SEIR2020,Harkin,Rame}.  In particular, if one knows the asymptotic behavior in the vicinity of the boundaries of the physical domain, the method of asymptotic approximants may be used to constrain analytic continuation by enforcing these behaviors at both ends, thus providing efficient and accurate analytic infinite series solutions which are in closed form when truncated. In cases where a power series solution converges over the whole domain of interest, but convergence is slow, the same approach may be used to re-sum the series \textit{exactly} and enhance convergence~\cite{Harkin, Rame}.  In this paper, we  establish that,  with a judicious choice of initial conditions, a convergent power series solution to the nonlinear pendulum ODE may be obtained.  The convergence of the series solution is then accelerated by exactly re-summing it in a form motivated by the asymptotic approximant technique.  We should note that the word `approximant' is not to be confused with `approximate' and is only used here to tie to previous literature.  The infinite series solution we present here, both in its original and re-summed form, is an exact solution.  

% the singularities of the pendulum were not known, one could still use knowledge of the asymptotic behavior at the two extreme ends of the domain to form an approximant that analytically continues the series.

While one goal of this work is to assemble the above-described solution elements from across the centuries in order to give a full treatment of the pendulum series for direct use, another aim here is to add this example to a growing list of problems~\cite{Barlow:2017,FS2020,SIR2020,SEIR2020,Harkin,Rame} in which power series solutions are demonstrated to be a viable solution technique for nonlinear ODEs, even in the more typical cases where singularity locations are unknown.  The paper is organized as follows. In section~\ref{sec:background} we detail the key structural elements of the pendulum equation, as well as its solution in terms of Jacobi elliptic functions. In section~\ref{sec:series}, we demonstrate an equivalent but alternative approach to~\cite{Fairen} for deriving the power series solution to the differential equation.  Using the known placement of the singularities in the complex time plane, we illustrate that the point of expansion \--- or equivalently, initial conditions \--- can be chosen such that the series solution encompasses the entire physical domain, using reflections about a judiciously chosen effective period. In section~\ref{sec:approximant}, the method of asymptotic approximants is used to motivate an exact resummation of the pendulum series, and details are provided on its construction and use. Appendix~\ref{sec:A1} provides an efficient way of implementing truncations of this resummation. In section~\ref{sec:period}, we develop an exact series resummation for the complete elliptic integral of the first kind that is used in period calculations.  This period resummation, in turn, is utilized to  periodically extend the solution forms of sections~\ref{sec:series} and~\ref{sec:approximant}. Concluding remarks are given in section~\ref{sec:conclusions} on the pendulum series solution provided here. This solution adds to the growing body of literature in which nonlinear ODEs are solved using series solutions, and so some closing remarks on this connection are provided as well.   

\section{Background \label{sec:background}}

\begin{figure}[h!]
    \centering
    \begin{tikzpicture}
    \coordinate (centro) at (0,0);
    \draw[dashed,gray,-] (centro) -- ++ (0,-3.5) node (mary) [black,below]{$ $};
    \draw[thick] (centro) -- ++(300:3) coordinate (bob);
    \draw (1.9, -2.6) node {$m$};
    % \draw[thick, ->] (1.5, -2.6) -- (1.5, -3.5);
    % \draw (1.9, -3.5) node {$mg$};
    \draw (1.1,-1.17) node {$\ell$};
    \fill (bob) circle (0.1);
    \pic [draw, ->, "$\theta$", angle eccentricity=1.5] {angle = mary--centro--bob};
    \end{tikzpicture}
    \caption{Setup of the simple pendulum problem. The mass, $m$, at the end of the string (of length $\ell$) swings freely from some initial angle and angular velocity. In this reference frame, the acceleration of gravity, $g$, acts vertically downward.}
    \label{fig:problem_setup}
\end{figure}
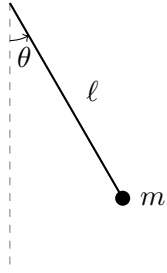
A schematic of the problem to be examined is provided in figure~\ref{fig:problem_setup}, which shows a pendulum bob of point mass $m$ anchored by a mass-less string of length $\ell$.  During motion, the bob is  displaced from the vertical by an angle $\theta$, acted upon by a downward gravitational force $mg$ ($g$ being the acceleration of gravity), and moves with an angular acceleration of $d^2\theta/dt^2$ where $t$ is time.  
Subject to the above assumptions, a moment balance leads to the ``simple pendulum equation''~\cite{Agnew}:
\begin{equation*}
    \frac{d^2\theta}{dt^2}+\frac{g}{\ell}\sin(\theta)=0.
    \label{pendeq}
\end{equation*} 
Defining the dimensionless time $\bar{t}=t\sqrt{\frac{g}{\ell}}$ as done in~\cite{Strogatz}, the initial value problem for the simple pendulum then becomes
\begin{subequations}
\begin{equation}
     \frac{d^2\theta}{d\bar{t}^2}+\sin(\theta)=0
     \label{pendalone}
\end{equation}
\begin{equation}
       \theta(0)=\theta_0,~~\frac{d\theta}{d\bar{t}}(0)=\bar{\omega}_0,
       \label{ics}
\end{equation}
    \label{pendeqnd}
\end{subequations}
where $\theta_0$ and $\bar{\omega}_0$ are the initial angular displacement and dimensionless angular velocity\footnote{The dimensionless angular velocity is defined in terms of its dimensional equivalent, $\omega_0$, as $\bar{\omega}_0=\omega_0\sqrt{\frac{\ell}{g}}$.}, respectively. An energy constant for the simple pendulum is obtained by multiplying both sides of~(\ref{pendalone}) by $d\theta/d\bar{t}$, integrating both sides with respect to $\bar{t}$, and applying~(\ref{ics}) to obtain
\begin{equation}
   \frac{1}{2}\left(\frac{d\theta}{d\bar{t}}\right)^2 + 1 - \cos(\theta)=\frac{1}{2}\bar{\omega}_0^2 + 1 - \cos(\theta_0)\equiv\bar{E}\ge0.
    \label{energynd}
\end{equation}
In~(\ref{energynd}), $\bar{E}$ is the dimensionless total energy\footnote{The dimensionless total energy is defined in terms of its dimensional equivalent, $E$, as $\bar{E}=\frac{E}{mg\ell}$.}, which is zero when the pendulum is resting at zero angle. 
% At any given point along the pendulum's trajectory, the sum of total energy, $E$, is constant and given by 
% \begin{equation}
% E = \frac{1}{2}m\left(\frac{d\theta}{dt}\right)^2\ell^2 + mg\ell[1 - \cos(\theta)],
% \label{energydim}
% \end{equation}
% where is $g$ is the acceleration of gravity and $t$ is time. 

%It is worth noting that analysis in this paper focuses primarily on the solution $\theta(t)$. 

% Prior studies more often focus on the inverse solution $t(\theta)$, as it can be obtained as an elliptic integral through application of separation of variables to~(\ref{energynd}).  However, the solution $t(\theta)$ is not typically used to describe the trajectory, but rather to provide an expression for the \textit{period} of the pendulum trajectory, using appropriate integration limits; see~\cite{Hinrichsen} for a recent and comprehensive review of these efforts.  

\begin{figure}[h!]
    \centering
    \includegraphics[width=5in]{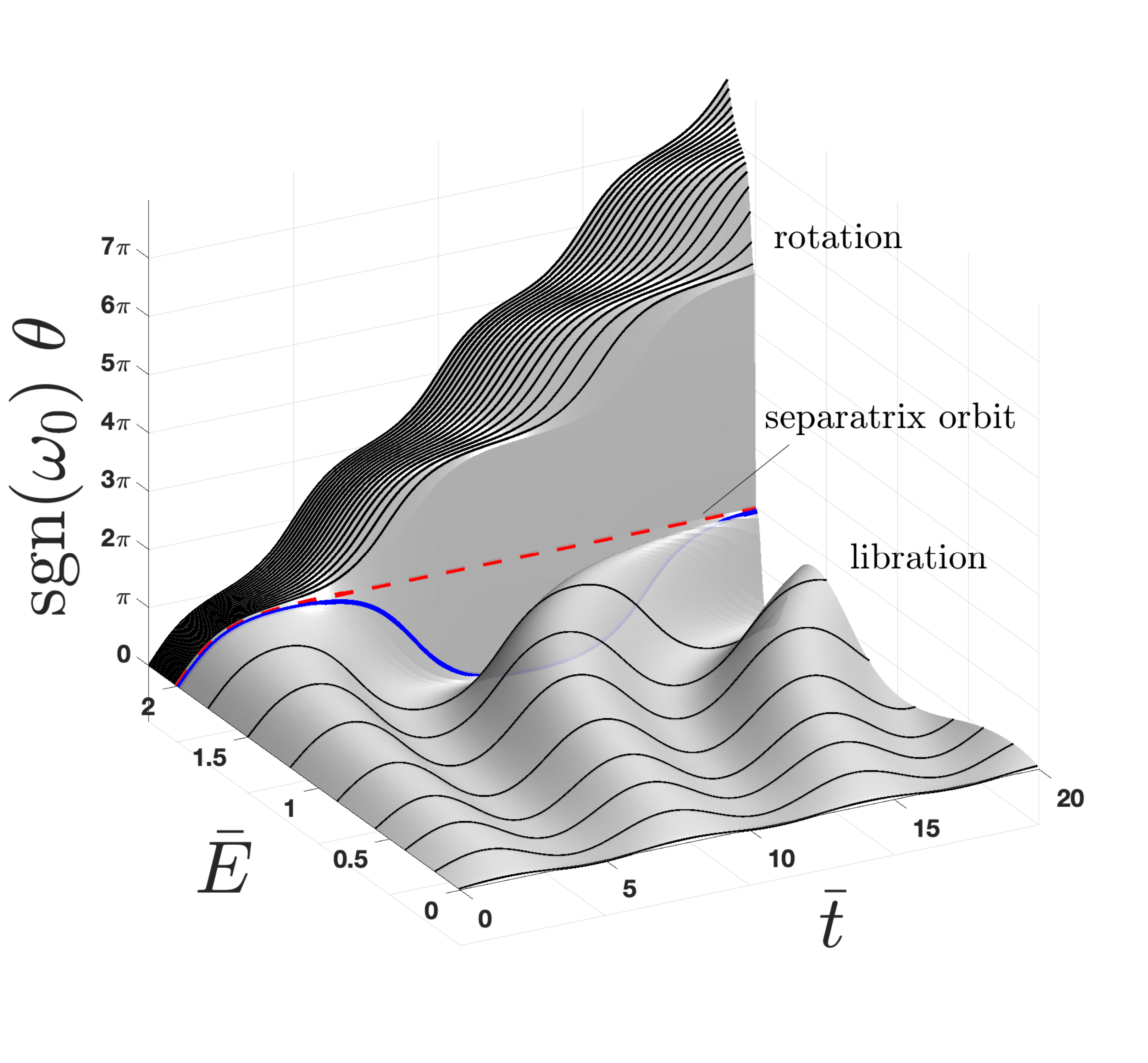}
    \caption{ The solution to~(\ref{pendeqnd}) describing the trajectory of a pendulum over time for several energy levels $\bar{E}$, encompassing the 3 possibilities: libration ($\bar{E}<2$), the separatrix orbit ($\bar{E}$=2), and rotation ($\bar{E}>2$).  The location of $\bar{t}=0$ is chosen to be at the bottom of the pendulum's trajectory, such that all energy states can be shown from the same starting point.    The  $\text{sgn}(\omega_0)$ is used as a pre-factor on $\theta$ to indicate that $\omega_0>0$ leads to the surface of solutions as shown while $\omega_0<0$ leads to a reflection of the surface about the $t-\bar{E}$ plane. The solutions shown may be considered exact to within machine precision and are generated using analytic formulae~(\ref{newseries}), ~(\ref{fullresummation}), and~(\ref{Knew}) developed in sections~\ref{sec:convergence},~\ref{sec:approximant}, and~\ref{sec:period}.}
    \label{fig:surf}
\end{figure}

% Continuing along the path of determining $\theta(t)$,
Equation~(\ref{pendeqnd}) \--- and the types of solutions it admits \--- have been studied extensively in prior literature~\cite{Whittaker,Lima,Ochs,chapter}.  All possible \textit{moving} (ignoring fixed points) solutions to~(\ref{pendeqnd}) are shown in figure~\ref{fig:surf}. Initial conditions used for figure~\ref{fig:surf} are purposely chosen to be $\theta_0$=0 while letting $\bar{\omega}_0$  vary to determine $\bar{E}$ according to~(\ref{energynd}) (i.e., pushing the pendulum from the bottom) since these constraints are consistent with all possible pendulum motion. For other combinations of initial conditions, one only needs to compute $\bar{E}$, locate the initial conditions in figure~\ref{fig:surf} (shifting the $\bar{t}$ origin to that point), and then follow the fixed $\bar{E}$ trajectory from that point forward \-- this is consistent with the properties of autonomous ODEs.   The  $\text{sgn}(\omega_0)$ in figure~\ref{fig:surf} is used as a pre-factor on $\theta$ to indicate that $\omega_0>0$ leads to the surface of solutions as shown while $\omega_0<0$ leads to a reflection of the surface about the $t-\bar{E}$ plane. 

The types of pendulum solutions admitted by~(\ref{pendeqnd}) and shown in figure~\ref{fig:surf} are classified by their value of $\bar{E}\in[0,\infty)$, approaching the standard closed-form linear approximations as $\bar{E}\to0$ and $\bar{E}\to\infty$.  As $\bar{E}\to0$, simple harmonic motion arises in the form $\theta\approx \theta_0\cos(\bar{t})+\bar{\omega}_0\sin(\bar{t})$, this limit being accessed for small $\theta$ in~(\ref{pendeqnd}) such that $\sin(\theta)\approx\theta$.   As $\bar{E}\to\infty$, uniform circular motion arises in the form $\theta\approx\bar{\omega}_0t$, this limit being accessed by making $d\theta/d\bar{t}$ large in~(\ref{energynd}) such that the potential energy term $1-\cos(\theta_0)$ (always of $O$(1)) can be neglected, allowing direct integration of~(\ref{energynd}). %While these limiting behaviors are unobtainable in the exact sense for any finite value of the solution, they are a point of familiarity for framing our discussion of the exact solution behaviors that follow. 

In reality, the system~(\ref{pendeqnd}) admits 3 types of pendulum motion, as shown in figure~\ref{fig:surf}: the \textit{librating} (i.e., ``swinging'' back and forth) trajectory for $\bar{E}<2$, the \textit{separatrix} orbit for $\bar{E}=2$ (dashed line in figure) in which the pendulum rises towards the top, never reaching it yet never falling back down, and the \textit{rotating} or ``overturning'' case for $\bar{E}>2$.

For definiteness, we define the pendulum's period $T$ as follows. For libration ($\bar{E}<2$), we use the usual definition of period, such that $\theta$ satisfies $\theta(\bar{t})=\theta(\bar{t}\pm T)$ while for rotation ($\bar{E}>2$) we use the unwound definition  $\theta(\bar{t})=\theta(\bar{t}\pm T)\mp \text{sgn}(\bar{\omega}_0)2\pi$ where $\text{sgn}(\bar{\omega}_0)$ is 1 for counter-clockwise rotation and $-1$ for clockwise rotation. For the separatrix orbit ($\bar{E}=2$), $T$ is infinite and the solution is not periodic. As mentioned in section~\ref{sec:intro}, the evaluation of $T$ itself dates back to the 18$^\mathrm{th}$ century~\cite{Stillwell,Gray}; it is obtained by applying separation of variables to~(\ref{energynd}), solving for $\bar{t}$ in terms of an integral in $\theta$, and then either: (for $\bar{E}<2$) letting $\bar{\omega}_0$=0, integrating between the quarter-period limits ($\theta=0$ to $\theta=\theta_0$), and multiplying by 4; or (for $\bar{E}>2$) letting $\theta_0=0$, integrating between the half-period limits ($\theta=0$ to $\theta=\pi$) and multiplying by 2. Although these integrals cannot be evaluated in closed-form, they may be manipulated into elliptic integral form by applying the half-angle identity, employing a variable substitution\footnote{For $\bar{E}<2$, let $\sin\frac{\theta_0}{2}\sin u=\sin\frac{\theta}{2}$. For $\bar{E}>2$, let $u=\frac{\theta}{2}$}, and a rewriting the resulting expression in terms of $\bar{E}$ to obtain 
\begin{equation}
T=
\begin{cases}
4~K\!\left(\sqrt{\frac{\bar{E}}{2}}\right) & \text{if $\bar{E}<2$ } \\
2\sqrt{\frac{2}{\bar{E}}}~K\!\left(\sqrt{\frac{2}{\bar{E}}}\right) & \text{if $\bar{E}>2$}.
  \end{cases}
  \label{T}
\end{equation}
where 
\begin{equation} K(k) = \int_0^\frac{\pi}{2} \frac{d\phi}{\sqrt{1 - k^2\sin^2(\phi)}}.
\label{K}
\end{equation}
In~(\ref{K}), $K(k)$ is the complete elliptic integral of the first kind~\cite{Byrd}.  

While figure~\ref{fig:surf} shows trajectory solutions of~(\ref{pendeqnd}) over several periods, all information above and below the separatrix for all $\bar{t}$ is, in fact, obtainable from reflections and translations of the solution in a finite interval $t\in[0,T^*]$, afforded by the symmetry of the problem. We define the \textit{effective period} $T^*$ to be the minimum portion of the period necessary to map out the entire solution. In the case of libration, the pendulum sweeps out the same (albeit negative) angles in time as it moves from  $\theta=0$ to $-\theta_0$ as it does going from $\theta=0$ to $\theta=\theta_0$. Thus $T^*$ for a librating pendulum is the quarter-period. In the case of rotation, instead of the back\--and\--forth motion described above, the pendulum swings ``over the top'' and ends up back where it started from. Here, the pendulum sweeps out the same angles (albeit with a $\pi$ offset) in time as it moves from $\theta=\theta_0$ to $\theta_0+~\text{sgn}(\bar{\omega}_0)\pi$ as it does going from $\theta=\theta_0 +\text{sgn}(\bar{\omega}_0)\pi$ to $\theta_0+2~\text{sgn}(\bar{\omega}_0)\pi$.  Thus $T^*$ for a rotating pendulum is the half-period.  To summarize:
\begin{equation}
T^*=
\begin{cases}
\frac{T}{4} & \text{if $\bar{E}<2$} \\
% \infty & \text{if $\bar{E}=2$} \\
\frac{T}{2} & \text{if $\bar{E}>2$},
  \end{cases}
  \label{Tstar}
\end{equation}
where the separatrix orbit ($\bar{E}=2$) is omitted from the above definition since it is not periodic.

Interestingly, initial conditions corresponding to $\bar{E}=2$ lead to simplifications such that the elliptic integral for $\bar{t}(\theta)$ (obtained by separating~(\ref{energynd})) can be evaluated in closed form and $\theta(t)$ can be solved for explicitly~\cite{Lima}\footnote{Note that, in arriving at~(\ref{critical}) through the route of elliptic integrals described here, $\theta_0$ is taken to be 0 in~\cite{Lima}, whereas we leave $\theta_0$ to be arbitrary. An alternative derivation (starting from elliptic \textit{functions}) is given in~\cite{Whittaker} where $\theta_0$ is also taken to be 0.} as
\begin{subequations}
\begin{equation}
\theta=-\pi+4~\text{arctan}\left[e^{\bar{t}}\tan\left(\frac{\theta_0+\pi}{4}\right)\right],\text{ for }~~\bar{E}=2.    
    \label{critical}
\end{equation}
For the libration and rotation cases, no known solution exists for~(\ref{pendeqnd}) in terms of elementary functions; however, as mentioned in section~\ref{sec:intro}, the solution can be written in terms of Jacobi elliptic functions.  These solutions are obtained by making the substitution $y = \sin(\theta/2)$ into the energy equation~(\ref{energynd}), enforcing initial conditions for either $\bar{E}<2$ or $\bar{E}>2$,  and recognizing the resulting ODE as that governing an elliptic function~\cite{chapter}.  For libration and rotation, the solutions are respectively, 
\begin{equation}
\theta = 2 \arcsin\left[\sqrt{\frac{\bar{E}}{2}}\operatorname{cd}\left(\bar{t}, \sqrt{\frac{\bar{E}}{2}}\right)\right],\text{ for }~~\bar{E}<2,~\bar{\omega}_0=0,   
\label{subcritical}    
\end{equation}
and 
\begin{equation}
\theta = 2 \arcsin\left[\operatorname{sn}\left(\sqrt{\frac{\bar{E}}{2}}\bar{t}, \sqrt{\frac{2}{\bar{E}}}\right)\right],\text{ for }~~\bar{E}>2,~\theta_0=0  
\label{supercritical}    
\end{equation}
\end{subequations}
for a given $\bar{E}$.  In~(\ref{critical}),~(\ref{subcritical}), and~(\ref{supercritical}), note that, once the energy $\bar{E}$ is set, only one initial condition (either $\theta_0$ or $\bar{\omega}_0$) is needed to specify the solution. The remaining initial condition may then be extracted from~(\ref{energynd}). As stated previously, all initial conditions can be accessed by shifting solutions of the same energy $\bar{E}$.  As there can be notational differences in the literature, for definiteness we utilize the convention~\cite{Byrd,chapter} in~(\ref{subcritical}) and ~(\ref{supercritical}) that the Jacobi elliptic functions $\operatorname{cd}(u,k)$ and $\operatorname{sn}(u,k)$ satisfy the nonlinear initial value problems 
\[\left(\frac{d\operatorname{cd}}{du}\right)^2=\left(1-\operatorname{cd}^2\right)\left(1-k^2\operatorname{cd}^2\right),~~\operatorname{cd}(0,k)=1\]
and
\[\left(\frac{d\operatorname{sn}}{du}\right)^2=\left(1-\operatorname{sn}^2\right)\left(1-k^2\operatorname{sn}^2\right),~~\operatorname{sn}(0,k)=0.\]
The solution of these equations is not known analytically in closed form; however, efficient numerical algorithms exist to evaluate the Jacobi elliptic functions in many scientific computing packages. As an alternative to the numerical solution of~(\ref{pendeqnd}) or the solutions forms~(\ref{subcritical}) and~(\ref{supercritical}), in section~\ref{sec:series} of this paper we generate an exact solution to~(\ref{pendeqnd}) via a convergent power series, whose convergence is subsequently accelerated in section~\ref{sec:approximant}.  Although the solution of the pendulum equation in this paper provide an alternative to the numerical evaluation of the elliptic functions, the expressions~~(\ref{subcritical}) and~(\ref{supercritical}) are important when making exact statements about singularity locations in section~\ref{sec:series}. Note that there are several other equivalent representations of the solution to~(\ref{pendeqnd}) involving other Jacobi elliptic functions.  In fact, one may also express the solution as a single elliptic function that governs all pendulum movement~\cite{Ochs}.

% Note that, although the libration and rotation regions in~\ref{fig:surf} could certainly have been generated by directly numerically solving~(\ref{pendeqnd}) or using~(\ref{subcritical}) and~(\ref{supercritical}) (elliptic functions are commonly included in scientific software and programming languages), we used our accelerated series formulation of the exact solution to~(\ref{pendeqnd}) to be described in section~\ref{sec:approximant} to (relatively quickly) generate the surface, as a testimony to the virtues of analytic series solutions. However, before discussing ways to accelerate convergence of the series, we need to discuss the series solution itself, given in the following section. 

\section{The Pendulum Series and Complex-Time Singularities \label{sec:series}}
The power series solution of~(\ref{pendeqnd}) follows a standard approach~\cite{Agnew}, and is developed in section~\ref{sec:reccurence}; this leads to a recurrence formula for its coefficients, equivalent to that obtained via a different technique in~\cite{Fairen}. However, the previous solution~\cite{Fairen} leaves implementation to the reader with arbitrary initial conditions and does not examine convergence properties of the series, leaving the solution itself as incomplete. The radius of convergence of the series based on closest singularities is determined in section~\ref{sec:sing}; through the placement of judicious expansion points (or, equivalently, initial conditions) this enables us to create a convergent series in the time range $t\in[0,T^*]$ (see~(\ref{Tstar}) and surrounding discussion). This solution is subsequently utilized to construct the solution to~(\ref{pendeqnd}) for all time in section~\ref{sec:convergence}.

\subsection{Power Series Solution to~(\ref{pendeqnd}) for $\theta(t)$ \label{sec:reccurence}}
The power series solution may be generated by assuming
\[\theta=\sum_{n=0}^\infty a_n \bar{t}^n,\]
substituting this into~(\ref{pendeqnd}), and equating like-powers of $\bar{t}$ on each side. The result is a recurrence expressed as
\begin{equation}
    a_{n+2} = \frac{- s_n}{(n+1)(n+2)},~a_0=\theta_0,~a_1=\bar{\omega}_0
    \label{pendcoeffs}
\end{equation}
where $s_n$ corresponds to coefficients one obtains by expanding the sine of an infinite series, i.e.,  
\begin{equation}
    \sin\left(\sum_{n=0}^\infty a_n \bar{t}^n\right)=\sum_{n=0}^\infty s_n \bar{t}^n.
    \label{sineseries}
\end{equation} 
The expansion~(\ref{sineseries}) is obtained using the same method as that for raising a a series to a power~\cite{Henrici}, taking the log of a series~\cite{SIR2020}, and taking the exponential of a series~\cite{SEIR2020}.  In fact, we make direct use of the latter formula in what follows, and repeat its construction here. The exponential of a series
\begin{equation}
    \text{exp}\left(\sum_{n=0}^\infty a_n \bar{t}^n\right)=\sum_{n=0}^\infty b_n \bar{t}^n
    \label{eseries}
\end{equation} 
is obtained by first noting that $b_0=e^{a_0}$ and then differentiating both sides of~(\ref{eseries}) such that 
\begin{equation}
\left[\sum_{n=0}^\infty (n+1)~a_{n+1}~\bar{t}^n\right]\text{exp}\left(\sum_{n=0}^\infty a_n \bar{t}^n\right)=\sum_{n=0}^\infty (n+1)~b_{n+1}~\bar{t}^n,
\end{equation}
which can be rewritten, using~(\ref{eseries}), as
\begin{equation}
    \left[\sum_{n=0}^\infty (n+1)~a_{n+1}~\bar{t}^n\right]\sum_{n=0}^\infty b_n \bar{t}^n=\sum_{n=0}^\infty (n+1)~b_{n+1}~\bar{t}^n.
    \label{deseries}
\end{equation}
Cauchy's product rule is then applied to the left-hand side of~(\ref{deseries}) to arrive at 
\begin{equation}
    \sum_{n=0}^\infty \left[\sum_{k=0}^n (k+1)~a_{k+1}~b_{n-k}\right] \bar{t}^n=\sum_{n=0}^\infty (n+1)~b_{n+1}~\bar{t}^n,
    \label{eCauchy}
\end{equation}
after which, equating like powers of $\bar{t}$ in~(\ref{eCauchy}) leads to a recursive expression for the coefficients of the expansion of an exponential of an infinite series, collected with~(\ref{eseries}) as
\begin{align}
\nonumber
    &\text{exp}\left(\sum_{n=0}^\infty a_n \bar{t}^n\right)=\sum_{n=0}^\infty b_n \bar{t}^n,\\ &b_{n+1}=\frac{1}{(n+1)}\displaystyle\sum_{k=0}^n (k+1)~a_{k+1}~b_{n-k},~~b_0=e^{a_0}.
    \label{eCoeffs}
\end{align}
To find the coefficients of~(\ref{sineseries}), we now define the cosine of a series as 
\begin{equation}
    \cos\left(\sum_{n=0}^\infty a_n \bar{t}^n\right)=\sum_{n=0}^\infty c_n \bar{t}^n,
    \label{cosseries}
\end{equation} 
 and apply Euler's identity
\[\text{exp}\left(i\sum_{n=0}^\infty a_n \bar{t}^n\right)=\cos\left(\sum_{n=0}^\infty a_n \bar{t}^n\right)+i\sin\left(\sum_{n=0}^\infty a_n \bar{t}^n\right),\]
directly to~(\ref{eCoeffs}) (using the notation of~(\ref{sineseries}) and~(\ref{cosseries})) to yield:
\begin{equation}
    c_{n+1}+i~s_{n+1}=\frac{1}{(n+1)}\displaystyle\sum_{k=0}^n (k+1)~i~a_{k+1}~\left(c_{n-k}+i~s_{n-k}\right).
    \label{components}
\end{equation}
In~(\ref{components}), it is clear that $s_0=\sin(a_0)$ and $c_0=\cos(a_0)$ from~(\ref{sineseries}) and~(\ref{cosseries}), respectively.  Before continuing, it is useful to note that $a_n$, $c_n$, and $s_n$ are real for all $n$. Now, equating the real part of each side of~(\ref{components}) we obtain a recursion for the coefficients of~(\ref{cosseries}). Similarly, equating the imaginary part of each side of~(\ref{components}), we obtain a recursion for the coefficients of~(\ref{sineseries}). Collecting these (coupled) recursions with~(\ref{pendcoeffs}), we may now write the series solution for the pendulum equation~(\ref{pendeqnd}) as
\begin{subequations}
\begin{equation}
    \theta=\sum_{n=0}^\infty a_n \bar{t}^n,~~|\bar{t}|<\bar{t}_\text{ROC}
    \label{sum}
\end{equation}
\begin{align}
    a_{n+2} &= \frac{- s_n}{(n+1)(n+2)},~a_0=\theta_0,~a_1=\bar{\omega}_0, \label{seriesrec} \\
        s_{n+1} &= \frac{1}{n+1}\sum_{k=0}^n (k+1)a_{k+1} c_{n-k},~s_0=\sin(\theta_0) \label{cosseriesrec}\\
    c_{n+1} &= \frac{-1}{n+1}\sum_{k=0}^n (k+1)a_{k+1} s_{n-k},~c_0=\cos(\theta_0). \label{sinseriesrec} 
\end{align}
 \label{pendseries}
 \end{subequations}
where $\bar{t}_\text{ROC}$ is the radius of convergence of the series, to be discussed in section~\ref{sec:sing}.  

An alternative derivation of the recursion in series~(\ref{pendseries}) is given in~\cite{Fairen} by writing~(\ref{pendeq}) as a coupled system of 4 ODEs in dependent variables $\theta$, $d\theta/dt$, $\sin(\theta)$, and $\cos(\theta)$ and relating the coefficients of their series expansions.  As another alternative to the approach used here to derive the series solution to~(\ref{pendeqnd}), one could substitute the well-known series expansions~\cite{Byrd,Wrigge} of the Jacobi elliptic functions into~(\ref{subcritical}) and~(\ref{supercritical}).  However, an arcsine of the resulting series is then required to extract the solution for $\theta(t)$, and for this reason, it is easier to proceed as presented above.

% As mentioned earlier, the reason for providing the details of the approach given by~(\ref{pendcoeffs}) through~(\ref{pendseries}) is that the identities may be used to develop series for modified pendula or other nonlinear ODEs.

% Also note that, instead of directly finding the series solution to~(\ref{pendeqnd}) as done above, one could utilize existing expansions~\cite{Byrd} for the elliptic functions and then substitute these into~(\ref{subcritical}) and~(\ref{supercritical}).  However, to arrive at the series for $\theta(\bar{t})$ from this approach, one would need to eventually take the expansion of the arcsine of the elliptic function series. This is possible, using similar series tricks as used above, but would arguably take just as much if not more time and paper as the direct approach that we take. Alternatively, if one used the expansion of the elliptic functions and took its arcsine without re-expanding, this would not provide a benefit in regards to the series divergence, as it would not remove the singularities (to be discussed in section~(\ref{sec:sing})). Besides, our work is meant to empower the use of series solutions nonlinear ODEs in general, not just those that can be represented by elliptic functions. 

\begin{figure}[h!]
\begin{tabular}{c}
(a)
\includegraphics[width=2.2in]{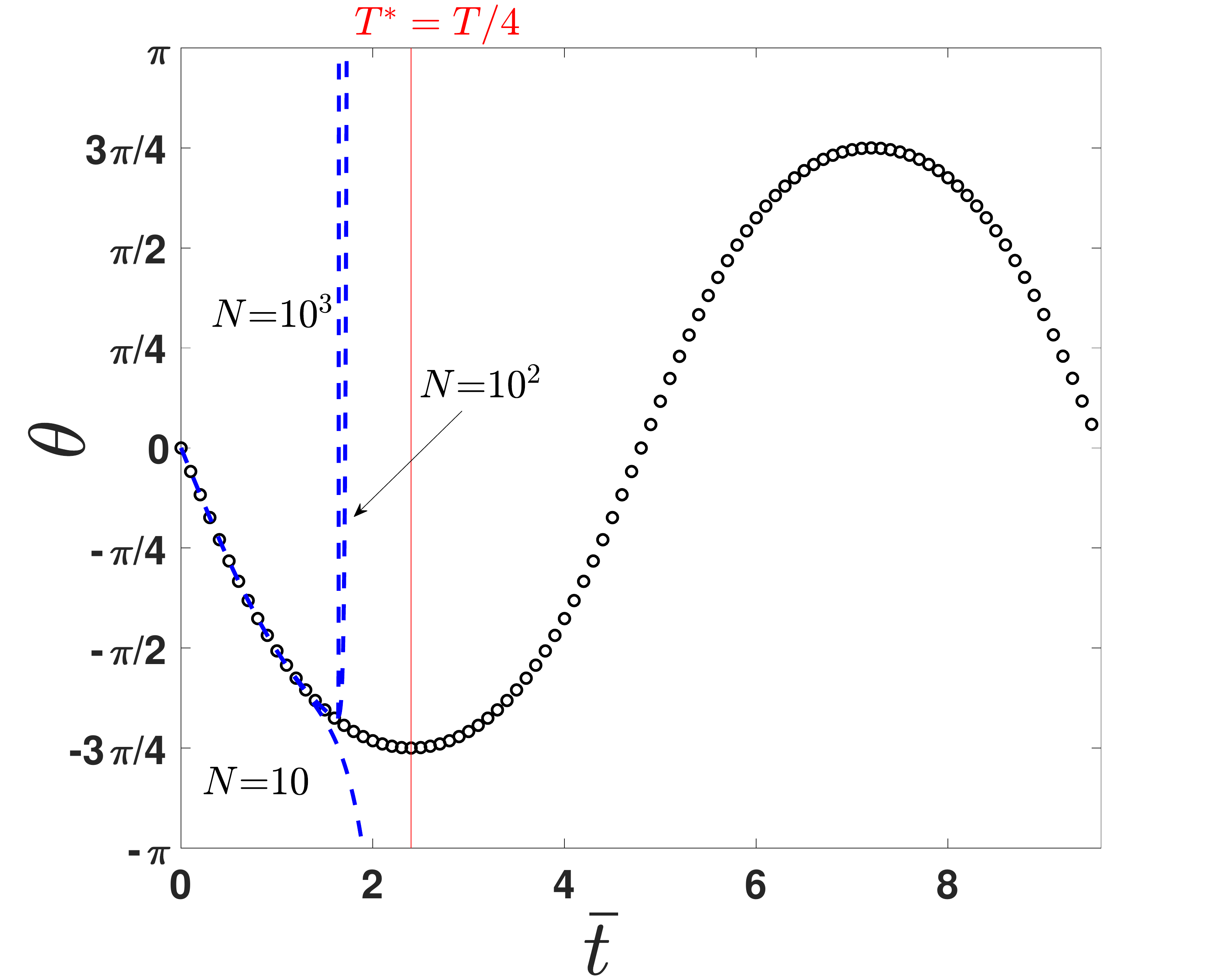}
(b)
\includegraphics[width=2.2in]{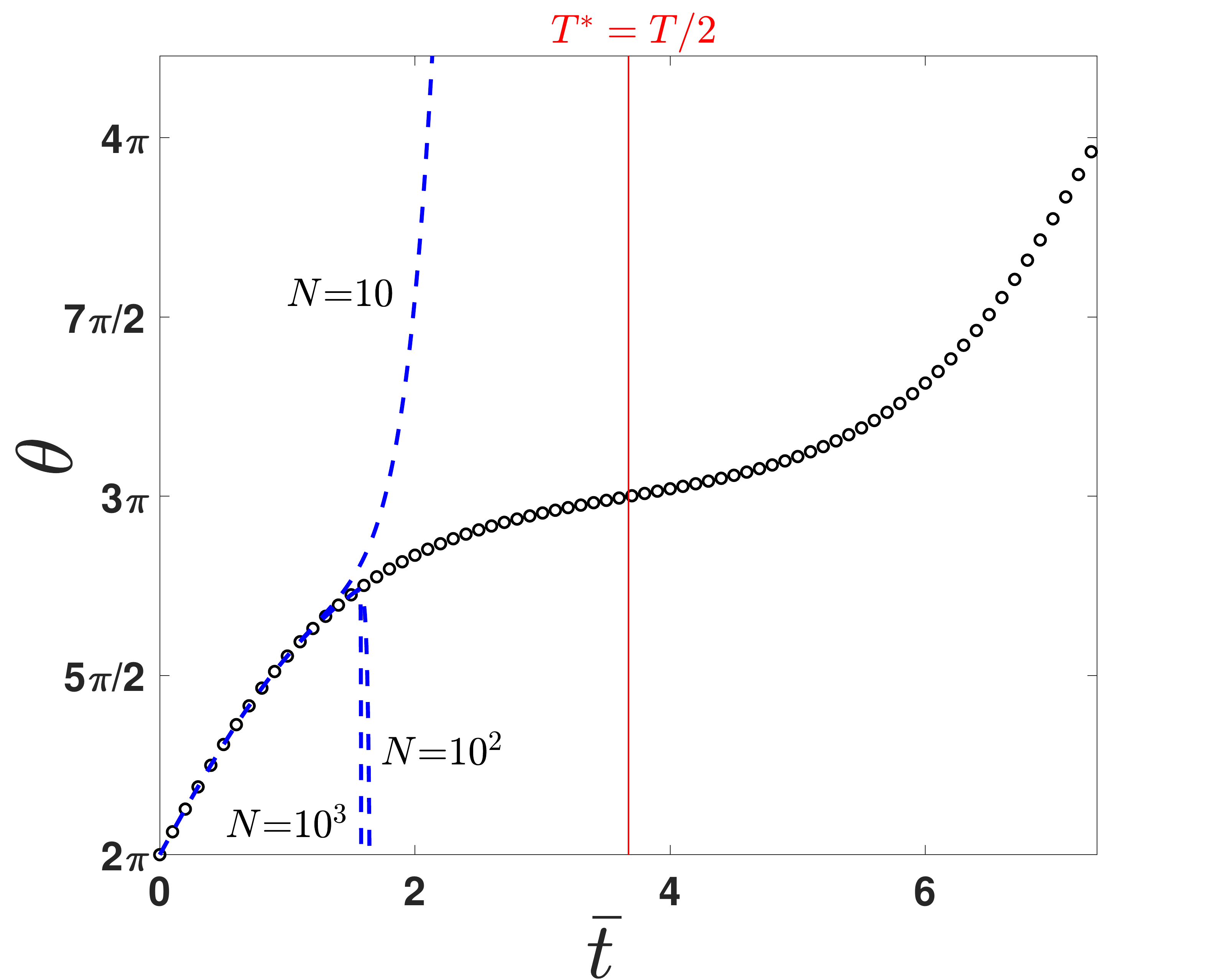}
\end{tabular}
    \caption{Series solution~(\ref{pendseries}) (dashed lines) taken up to (and including) the $\bar{t}^N$ term, shown for $N$=10,~10$^{2}$, and 10$^{3}$ and compared with the numerical solution ($\circ$) to~(\ref{pendeqnd}) over 1 period. The solution is initialized at the bottom of the pendulum's trajectory for cases of (a) libration ($\bar{E}<2$) with $\bar{E}=1.71$; and (b) rotation ($\bar{E}>2$) with $\bar{E}=2.02$. Note that, for either case, the series does not converge to the effective period $\bar{t}=T^*$ (indicated by a vertical line in the plots). }
    \label{fig:series}
\end{figure}

A comparison of the numerical solution\footnote{All numerical solutions in this work are generated using the 4$^\mathrm{th}$-order Runge-Kutta method with $\Delta\bar{t}=10^{-5}$. Points on figures have been down-sampled for readability.} to~(\ref{pendeqnd}) and the $(N+1)$-term partial sums of the series solution~(\ref{pendseries}) are shown in figure~\ref{fig:series} for both libration and rotation.  Note that the series diverges at a finite $\bar{t}_\text{ROC}$.   In generating the expansion~(\ref{pendseries}) for figure~\ref{fig:series}, we chose initial conditions at the bottom of the pendulum trajectory (i.e., $\theta_0=0$ in figure~\ref{fig:series}a). Figure~\ref{fig:series} indicates that for these choices of initial conditions, the pendulum series~(\ref{pendseries}) cannot fully describe the solution over the effective period $T^*$ due to divergence.  In the following section, we examine how $\bar{t}_\text{ROC}$ is affected by initial conditions, and determine an optimal choice to ensure convergence over the full range 0 to $T^*$.  Note that such an exploration is not restrictive, since the pendulum trajectory is invariant for a given choice of $\bar{E}$ as discussed previously \-- the locus of all initial conditions that achieve a given $\bar{E}$ yield trajectories that are simply shifted in time but sweep through the same angles.

% The case of libration in figure~\ref{fig:series}a, the  $\bar{E}$ $(\theta,\bar{\omega})=(3\pi/4,0)$ occurs, but is instead initialized at the bottom of this trajectory where ($\theta_0,\bar{\omega}_0$)=$\left(0,-\sqrt{2-\sqrt{2}}\right)$. The rotating case in figure~\ref{fig:series}b is also initialized at the bottom of the swing, one rotation around at ($\theta_0,\bar{\omega}_0$)=(2$\pi$, 2.01).  These choices are intentional and the value of $\bar{t}_\text{ROC}$ is set by singularities in the complex-$\bar{t}$ domain whose distance from the expansion point is a function of these initial conditions.  We now examine these singularities in order to overcome the apparent convergence barrier for all librating and rotating solutions of the pendulum problem. 

\subsection{\label{sec:sing}Singularities of $\theta(t)$ and series convergence}
Although any initial conditions $(\theta_0,\omega_0)$ may be chosen in~(\ref{pendseries}), the value of $\bar{t}_\text{ROC}$ is affected by this choice.  Here,  we show that there is an optimal choice of $(\theta_0,\omega_0)$ such that~(\ref{pendseries}) converges beyond the effective period $(\bar{t}=T^*)$ for all energy levels $\bar{E}$ \--- and thus all possible motions \--- of the pendulum.  This enables the series solution~(\ref{pendseries}) to be used as an exact solution that describes the pendulum trajectory for all $\bar{t}$.

In order to determine the radius of convergence, $\bar{t}_\text{ROC}$, of the pendulum series~(\ref{pendseries}), the singularities of the pendulum function in the complex-time plane are examined. The value of $\bar{t}_\text{ROC}$ is equal to the distance between the expansion point, chosen to be $\bar{t}$=0 in~(\ref{pendseries}), and its closest singularity in the complex plane. For the case of the separatrix orbit ($\bar{E}$=2), the arctan in~(\ref{critical}) admits logarithmic branch point singularities, the smallest of which are at $\bar{t}=\pm i\pi/2$ (for $\theta_0=0$ in~(\ref{critical})), such that $\bar{t}_\text{ROC}=\pi/2$.  However, since an exact solution is known for the separatrix case in terms of elementary functions, this restriction is illustrative but not relevant. Unlike the separatrix case for which a closed-form solution exists, or for the inverse solution $t(\theta)$ where the locations of singularities are known from the integrands in the elliptic integrals involved (such as~(\ref{K})), the time singularities of $\theta(t)$ for $\bar{E}\neq2$ are not obvious.  As mentioned in section~\ref{sec:intro}, it is typical that nonlinear ODEs such as~(\ref{pendeqnd}) contain ``spontaneous'' singularities that are not discernible by inspection of the ODE itself~\cite{bender,Barlow:2017,FS2020,SIR2020,SEIR2020}.  In such problems, an estimate for $\bar{t}_\text{ROC}$ may be obtained from numerical versions of the root or ratio tests, the latter of which is often represented graphically as a Domb-Sykes plot~\cite{DombSykes}.  

\begin{figure}[h!]
    % \centering
    \includegraphics[width=5in]{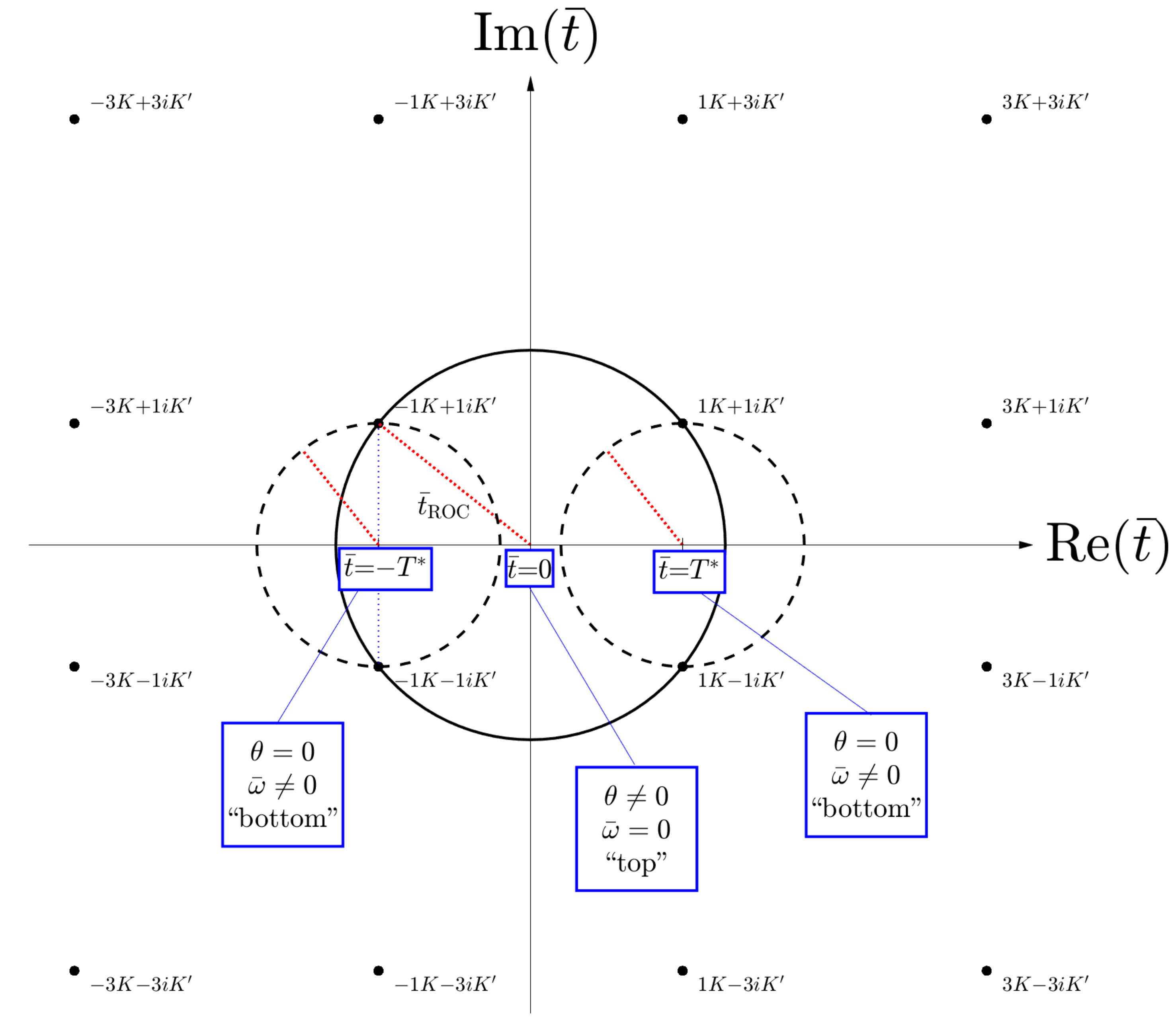}
    \caption{Schematic showing periodic structure of the singularities (shown as $\bullet$'s) of~(\ref{subcritical}) for the librating pendulum at a constant value of $\bar{E}<2$, as specified by~(\ref{libratingpoles}). The integrals $K(k)$ and $K'(k)$ are defined by~(\ref{K}) and~(\ref{Kp}), respectively.  The radius~(\ref{librationROC}) and placement of the solid circle corresponds to the original initial conditions given by~(\ref{subcritical}), initializing the pendulum at the top of its trajectory with zero speed; this is  the optimal radius of convergence $\bar{t}_\text{ROC}$ for the series solution~(\ref{pendseries}) of a librating pendulum.  By contrast, dashed circles show the radius of convergence when a librating pendulum is initialized at the bottom of the trajectory with nonzero angular speed.  }
    \label{fig:subcrit}
\end{figure}

Since the Jacobi elliptic functions have been well studied, their pole locations are known exactly, and thus the singularities in the pendulum function are also known as logarithmic branch points\footnote{The change of singularity type from a pole to a logarithmic branch point is due to the $\arcsin$ in ~(\ref{subcritical}) and~(\ref{supercritical}).} at the same locations through~(\ref{subcritical}) and~(\ref{supercritical}). This enables an exact determination of $\bar{t}_\text{ROC}$ for the pendulum series~(\ref{pendseries}).  For libration, the Jacobi elliptic function $\operatorname{cd}(\bar{t},\sqrt{\bar{E}/2})$ appearing in~(\ref{subcritical}) has an infinite number of poles~\cite{Byrd} at the complex time values 
\begin{equation}
\bar{t}_\mathrm{poles}=n_o \underbrace{K\left(\sqrt{\frac{\bar{E}}{2}}\right)}_{T^*=T/4} + n_o'~i K'\left(\sqrt{\frac{\bar{E}}{2}}\right),\text{ for }~~\bar{E}<2,~\bar{\omega}_0=0
\label{libratingpoles}
\end{equation}
where $n_o$ and $n_o'$ are odd integers, $K(k)$ is given by~(\ref{K}), and $K'(k)$ is (in general) defined as
\begin{equation}
    K'(k)\equiv K(\sqrt{1-k^2}).
    \label{Kp}
\end{equation}
In~(\ref{Kp}), $K'(k)$ is the magnitude of the ``imaginary quarter-period''~\cite{chapter,Appell}. The singularities prescribed by~(\ref{libratingpoles}) are shown in the complex $\bar{t}$-plane in figure~\ref{fig:subcrit} for a constant value of $\bar{E}<2$; note that the singularity structure relative to $\bar{t}=0$  is qualitatively the same for all values of $\bar{E}<2$ so one figure suffices to draw general conclusions. The series expansion about $t=0$ given by~(\ref{pendseries}) encounters the closest singularity when $n_o=n_o'=1$ in~(\ref{libratingpoles}) and the distance to it (i.e., the radius of convergence), $\bar{t}_\text{ROC}$, is given as follows:
\begin{equation}
    \bar{t}_\text{ROC}=\sqrt{{T^*}^2+\left[K'\left(\sqrt{\frac{\bar{E}}{2}}\right)\right]^2}> T^*\text{; valid for }~~\bar{E}<2,~\bar{\omega}_0=0.
    \label{librationROC}
\end{equation}
As noted in~(\ref{librationROC}), the radius of convergence of~(\ref{pendseries}) is always larger than the effective period $T^*$ (defined in~(\ref{Tstar})) for the given initial conditions, since the integral $K'(k)$ is never 0.  This circle of convergence is shown as a solid curve in figure~\ref{fig:subcrit} centered around the expansion point $\bar{t}$=0 with initial conditions at the \textit{top} of the pendulum's trajectory ($\theta_0\neq0,\bar{\omega}_0=0$); this is an optimal expansion point with respect to convergence for a given $\bar{E}<2$, as indicated in the figure.   If one were to expand about $\bar{t}=\pm T^*$ or, equivalently, redefine $\bar{t}=0$ on the real-axis and choose initial conditions at the \textit{bottom} of the trajectory ($\theta_0=0,\bar{\omega}_0\neq0$) for the same $\bar{E}<2$, the radius of convergence of~(\ref{pendseries}) is less than the effective period $T^*$, as shown by the smaller dashed circles in figure~\ref{fig:subcrit}; this is observed in~ figure~\ref{fig:series}a, where the series diverges prior to reaching $t=T^*$.  Note that the radius of convergence decreases monotonically for a given $E<2$ as the initial conditions are moved from the top of the trajectory ($\theta_0\neq0,\bar{\omega}_0=0$) to the bottom of the trajectory  ($\theta_0=0,\bar{\omega}_0\neq0$)  for the case of a librating pendulum.  This statement is made in the context of Equation~(\ref{energynd}), where $\theta_0$ and $\bar{\omega}_0$ are not independent quantities for a given $\bar{E}$ \--- hence the omission of $\theta_0$ in~(\ref{librationROC}).

For the case of rotation ($E>2$), a similar expression for pole locations is available~\cite{Byrd} for the Jacobi elliptic function $\operatorname{sn}(\sqrt{\bar{E}/2}~\bar{t},\sqrt{2/\bar{E}})$ that appears in~(\ref{supercritical}); here, the infinite number of singularities in the pendulum function are placed in the complex $\bar{t}$-plane at
\begin{equation}
    \bar{t}_\text{poles}=n_e\underbrace{\tilde{K}\left(\sqrt{\frac{2}{\bar{E}}}\right)}_{T^*=T/2}+n_o'~i\tilde{K}'\left(\sqrt{\frac{2}{\bar{E}}}\right),\text{ for }~~\bar{E}>2,~\theta_0=0,
    \label{rotatingpoles}
\end{equation}
where $n_e$ and $n_o'$ respectively denote even and odd integers and we define  $\tilde{K}(k)\equiv(\sqrt{2/\bar{E}})~K(k)$ and $\tilde{K'}(k)\equiv(\sqrt{2/\bar{E}})~K'(k)$ with $K'(k)$ given by~(\ref{Kp}).  The singularities prescribed by~(\ref{rotatingpoles}) are shown in the complex $\bar{t}$-plane in figure~\ref{fig:supercrit} for a constant value of $\bar{E}>2$.  The singularity structure for rotation is qualitatively similar for all values of $\bar{E}>2$ so one figure suffices to draw general conclusions.  The circle of convergence for the series expansion~(\ref{pendseries}) taken about $t=0$, with initial conditions specified by~(\ref{rotatingpoles}) at the bottom of the pendulum trajectory, is shown as the leftmost dashed circle in figure~\ref{fig:supercrit}.  A circle of convergence of the same size is shown centered around $t=T^*$, where the pendulum is also at the bottom of its trajectory.  In both cases,  the radius of convergence\footnote{The precise value of the radius of the smaller circles in figure~\ref{fig:supercrit} may be obtained by letting $n_e=0$ and $n_0'=1$ in~(\ref{rotatingpoles}) and taking its modulus; this radius is only larger than $T^*$ if $\bar{E}>4$. Thus, the pendulum series for rotation cannot converge beyond $\bar{t}=T^*$ for $\bar{E}\in[2,4]$ if initialized at the bottom of the trajectory.} is shown to be less than $T^*$, such that the series diverges prior to reaching the effective period. This is observed in~ figure~\ref{fig:series}b, where the series~(\ref{pendseries}) is initialized with $\theta_0=2\pi$ and diverges prior to reaching $t=T^*$.  All other points initialized at the bottom of a rotating pendulum's trajectory are equally detrimental to convergence for a given $\bar{E}>2$.  As shown in figure~{\ref{fig:supercrit}, an optimal expansion point (or initial condition) for a given energy $\bar{E}>2$ with respect to convergence is one that is again taken at the \textit{top} of the pendulum's trajectory (e.g., $\theta_0=\pi,\bar{\omega}_0\neq0$). Using the solid circle centered at $t=T^*$ and the surrounding lattice in figure~{\ref{fig:supercrit}, the optimal radius of convergence and initial conditions for rotation ($\bar{E}>2$) is given by
\begin{equation}
    \bar{t}_\text{ROC}=\sqrt{{T^*}^2+\left[\tilde{K}'\left(\sqrt{\frac{2}{\bar{E}}}\right)\right]^2}> T^*\text{; valid for }~~\bar{E}>2,~\theta_0=\pi,
    \label{rotationROCgood}
\end{equation}
where again, the radius of convergence is always larger than the effective period $T^*$, since $\tilde{K}'(k)$ is nonzero for finite $\bar{E}$. As a reminder when interpreting these equations, note that, by inspection of~(\ref{energynd}), $\theta_0$ and $\bar{\omega}_0$ are not independent quantities for a given $\bar{E}$, allowing for the omission of $\bar{\omega}_0$ in~(\ref{rotationROCgood}).

% The radius of convergence of~(\ref{pendseries}) with $\bar{E}>2$ and initial conditions corresponding to~(\ref{supercritical}) is thus
% \begin{equation}
%     \bar{t}_\text{ROC}=\tilde{K}'\left(\frac{2}{\bar{\omega}_0}\right)=\frac{2}{\bar{\omega}_0}{K}\left[\sqrt{1-\left(\frac{2}{\bar{\omega}_0}\right)^2}\right]\text{; valid for }~~\bar{E}>2,~\theta_0=0.
%     \label{rotationROC}
% \end{equation}
%  Equation~(\ref{rotationROC}) shows that, when using the initial condition $\theta_0=0$ for $\bar{E}>2$, the radius of convergence of~(\ref{pendseries}) is only larger than $T^*$ if $|\bar{\omega}_0|>2\sqrt{2}$.  This circle of convergence is shown in figure~{\ref{fig:supercrit} centered around the expansion point $\bar{t}$=0 with initial conditions at the \textit{bottom} of the pendulum's trajectory ($\theta_0=0,\bar{\omega}_0\neq0$); this point, along any other points initialized at the bottom of the trajectory (e.g., $\theta_0=2\pi$), are equally the worst expansion points with respect to convergence for a given $\bar{E}>2$, as indicated and deduced from the lattice shown. This is observed in~ figure~\ref{fig:series}b, where the series~(\ref{pendseries}) is initialized with $\theta_0=2\pi$ and diverges prior to reaching $t=T^*$. 

%Also, a substitution has been done in~(\ref{rotationROCgood}) such that the argument of $\tilde{K}'$ is written in terms of the new $\bar{\omega}_0$ at the shifted origin. 

\begin{figure}[h!]
    % \centering
    \includegraphics[width=5in]{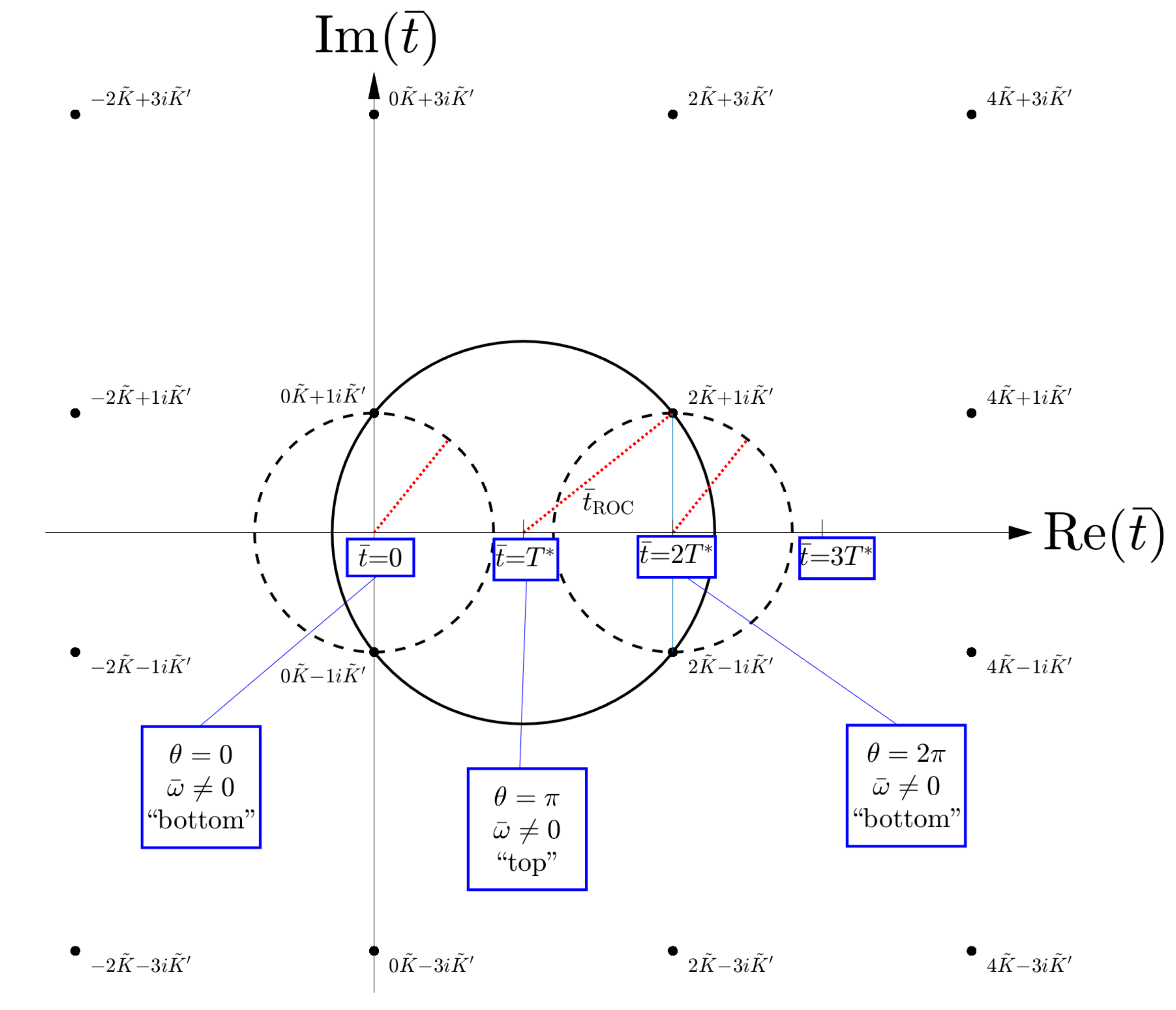}
    \caption{Schematic showing the periodic structure of the singularities (shown as $\bullet$'s) of~(\ref{subcritical}) of~(\ref{supercritical}) for the rotating pendulum at a constant value of $\bar{E}>2$, as specified by~(\ref{rotatingpoles}). $\tilde{K}(k)\equiv\sqrt{2/\bar{E}}~K(k)$ and $\tilde{K'}(k)\equiv\sqrt{2/\bar{E}}~K'(k)$, where the integrals $K(k)$ and $K'(k)$ are defined by~(\ref{K}) and~(\ref{Kp}), respectively. The leftmost dashed circle corresponds to the original initial conditions given by~(\ref{supercritical}), initializing the pendulum at the bottom of its trajectory.  The radius~(\ref{rotationROCgood}) and placement of the solid circle correspond to a shifted initialization at the top of the pendulum's trajectory, which is the optimal radius of convergence $\bar{t}_\text{ROC}$ for the series solution~(\ref{pendseries}) of a rotating pendulum.}
    \label{fig:supercrit}
\end{figure}

In both the case of libration and rotation, initializing the series solution to the pendulum at the highest point along the trajectory leads to the largest radius of convergence.  This allows us to exploit the symmetry of the problem such that the series solution~(\ref{pendseries}), when prescribed initial conditions as follows, 
\begin{equation}
    \theta(\bar{t})=\sum_{n=0}^\infty a_n t^n,~(\theta_0,~\bar{\omega}_0)=\begin{cases} \begin{rcases}
\left(\arccos(1-\bar{E}),~0\right), & \text{$\bar{E}<2$} \\
\left(\pi,~\pm\sqrt{2\bar{E}-4}\right),  & \text{$\bar{E}>2$},
 \end{rcases}, \end{cases},|\bar{t}|\le T^\star
    \label{seriesfixed}
\end{equation}
is an exact representation of the pendulum trajectory for all $t$ between 0 and $T^*$ and for any $\bar{E}$ describing librating and rotating motions of the pendulum.  In~(\ref{seriesfixed}), the $a_n$ coefficients are given in~(\ref{pendseries}), the $\pm$ symbol specifies either counterclockwise or clockwise rotation, and~(\ref{energynd}) has been used to represent the initial conditions (taken at the top of the trajectory) in terms of $\bar{E}$.

%Thus, if a series expansion for $\theta(t)$ about $t=0$ is obtained, its radius of convergence is given by~(\ref{eq:roc}) and will always be greater that the quarter-period. \TODO{Make the proclamations about exactness of the series}.  Such a series has been obtained in~\cite{Fairen}, although labeled there as being approximate. The series is given below, albeit written in a more compact, having been obtained by a slightly different (but mathematically equivalent) route.    

\subsection{Periodic extension of the series solution \label{sec:convergence}}
The convergent series solution~(\ref{seriesfixed}), justified in section~\ref{sec:sing} for $\bar{t}\le T^*$, allows for a convergent solution to~(\ref{pendeqnd}) for all $\bar{t}$ as follows.  For any given initial condition $(\theta_0,\bar{\omega}_0$), one first computes the energy $\bar{E}$ given by~(\ref{energynd}).  Using this value of $\bar{E}$, the convergent series for the solution up to the effective period $T^*$ is constructed from~(\ref{seriesfixed}) (or the resummation~(\ref{fullresummation}) provided in section~\ref{sec:approximant}) as 
\begin{subequations}
\begin{equation}
    \tilde{\theta}(\bar{t})\equiv\theta(\bar{t}\le T^*)
    \label{Ttilde}.
\end{equation}
The solution may then be extended for all $\bar{t}$ using the symmetry of the problem described in section~\ref{sec:background}, such that the periodic solution is given by
\begin{align}
\nonumber
  \theta(\bar{t})=&
\begin{cases}
  \begin{rcases}
  \tilde{\theta}(\hat{t}),& 0\le\hat{t}\le T^*\\
  -\tilde{\theta}(2T^*-\hat{t}),& T^*\le\hat{t}\le 2T^*\\
  -\tilde{\theta}(\hat{t}-2T^*),& 2T^*\le\hat{t}\le 3T^*\\
   \tilde{\theta}(4T^*-\hat{t}),& 3T^*\le\hat{t}\le 4T^*   
  \end{rcases},& \bar{E}<2\\
  ~\\
  \begin{rcases}
 \pm\frac{2\pi}{T}(\bar{t}-\hat{t})+\tilde{\theta}(\hat{t}),& 0\le\hat{t}\le T^*\\
\pm\frac{2\pi}{T}(\bar{t}-\hat{t}) -\tilde{\theta}(2T^*-\hat{t}) + 2\pi\pm2\pi,& T^*\le\hat{t}\le 2T^*\\
  \end{rcases}, & \bar{E}>2
\end{cases}\\ \nonumber~\\&\text{  where  } \hat{t}=t\text{ modulo } T,~~~~\bar{t}\ge0.
    \label{symmetry}
\end{align}
\label{newseries}
\end{subequations}
% \begin{align}
% \nonumber
%   \theta(\bar{t})=&
% \begin{cases}
%   \begin{rcases}
%   \tilde{\theta}(\hat{t}),& 0\le\hat{t}\le T^*\\
%   -\tilde{\theta}(2T^*-\hat{t}),& T^*\le\hat{t}\le 2T^*\\
%   -\tilde{\theta}(\hat{t}-2T^*),& 2T^*\le\hat{t}\le 3T^*\\
%   \tilde{\theta}(4T^*-\hat{t}),& 3T^*\le\hat{t}\le 4T^*   
%   \end{rcases},& \bar{E}<2\\
%   ~\\
%   \begin{rcases}
%  \frac{2\pi}{T}\text{sgn}(\bar{\omega}_0)(\bar{t}-\hat{t})+\tilde{\theta}(\hat{t}),& 0\le\hat{t}\le T^*\\
% \frac{2\pi}{T}\text{sgn}(\bar{\omega}_0)(\bar{t}-\hat{t}) -\tilde{\theta}(2T^*-\hat{t}) + 4\pi \delta_{\operatorname{sgn}(\bar{\omega}_0)}^1,& T^*\le\hat{t}\le 2T^*\\
%   \end{rcases}, & \bar{E}>2
% \end{cases}\\ \nonumber~\\&\text{  where  } \hat{t}=\operatorname{mod}(\bar{t}, T),~~~~\delta_{\operatorname{sgn}(\bar{\omega}_0)}^1=\begin{cases} 1, & \bar{\omega}_0>0\\ 0, & \bar{\omega}_0<0 \end{cases},~~~~\bar{t}\ge0
%     \label{symmetry}
% \end{align}
% \label{newseries}
% \end{subequations}
As in~(\ref{seriesfixed}), the $\pm$ symbol in~(\ref{symmetry}) specifies either counterclockwise or clockwise rotation.  The modulo operation\footnote{The modulo operation is available in most programming languages.  For example, the MATLAB syntax is mod$(\bar{t},T)$.} in~(\ref{symmetry}) is used such that $\hat{t}$ may be interpreted as the additional time beyond the smallest integer multiple of the period $T$. In section~\ref{sec:period}, we provide an exact and rapidly converging series for the computation of $T$ (and by extension, $T^*$), thus preserving the analyticity of the solution when used in~(\ref{newseries}).  The re-initializing and period extension given by~(\ref{newseries}) is shown in figures~\ref{fig:librationfixed} through~\ref{fig:librationseries}. Note that the formerly-divergent series solutions of figures~\ref{fig:series}a and~\ref{fig:series}b have been re-initialized using~(\ref{seriesfixed}) such that they are now convergent, and whose prior domains are now contained within the periodically extended domains of figures~\ref{fig:librationfixed} and~\ref{fig:rotationfixed}.  A final step in this solution process is to shift the origin such that the originally intended initial conditions are at $\bar{t}=0$. 

The norm of the error for $\bar{t}\in[0,T^*]$ incurred by series~(\ref{seriesfixed}) is plotted versus $\bar{E}$ (shown as dashed lines) in figure~\ref{fig:error} for various values of series truncation.   Although the error increases as $\bar{E}\to2$, for any given $\bar{E}$, the error can be reduced to machine precision by increasing the number of terms in the series.  This is expected, since there are no singularities of the pendulum function in the circle $|\bar{t}|<T^*$ for the initial conditions used in~(\ref{seriesfixed}).  
\begin{figure}[h!]
    \includegraphics[width=5.1in]{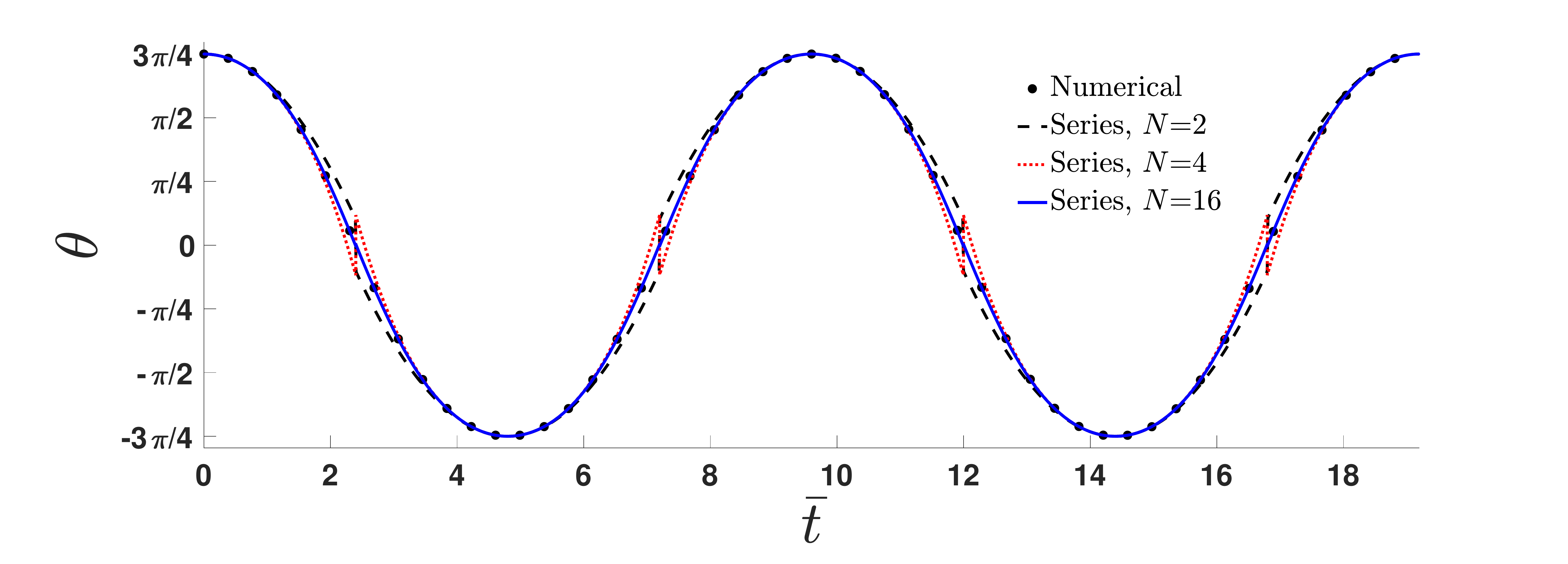}
    \caption{Series solution~(\ref{newseries}) (using~(\ref{seriesfixed}) for $\tilde{\theta}(t)$) taken up to (and including) the $\bar{t}^N$ term, shown for $N$=2, 4, and 6 and compared with the numerical solution ($\bullet$) to~(\ref{pendeqnd}) for the case of  libration ($\bar{E}<2$). The energy $\bar{E}=1.71$ is same as that in figure~\ref{fig:series}a but here the series is initialized at the top of the trajectory such that it converges. The solution is shown over 2 periods above, encompassing the solution in~figure~\ref{fig:series}a. }
    \label{fig:librationfixed}
\end{figure}

\begin{figure}[h!]
    \includegraphics[width=5.1in]{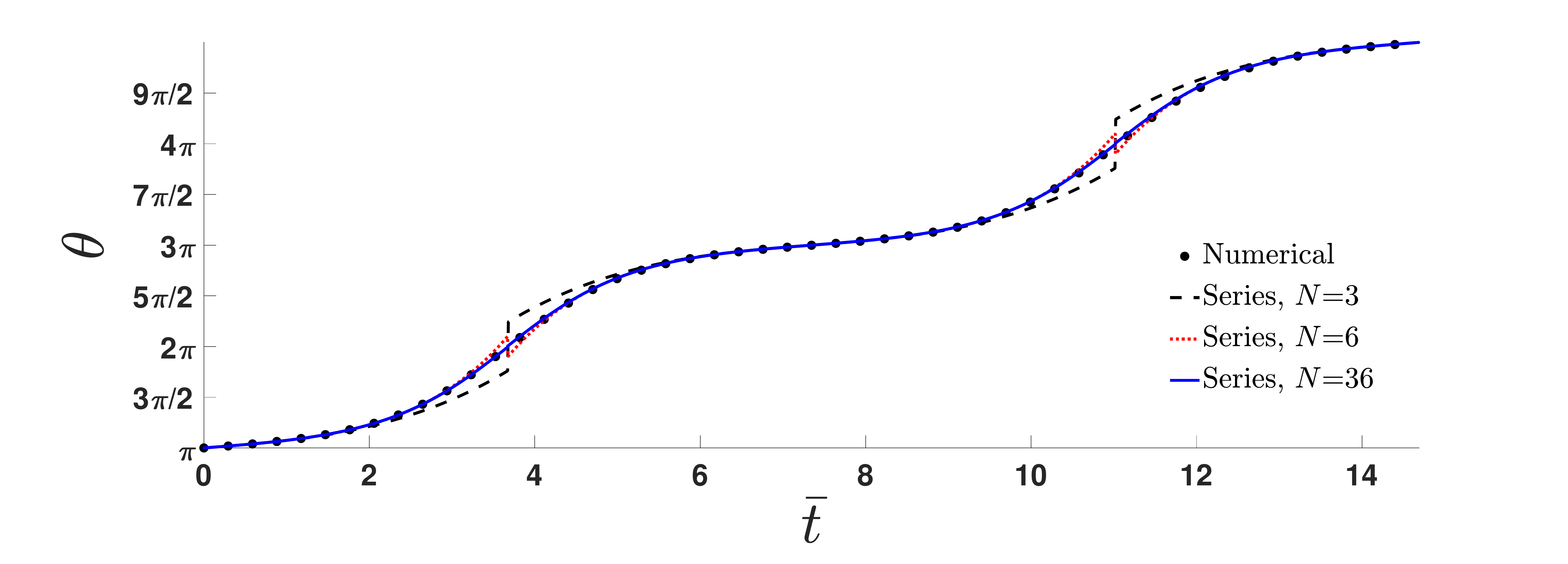}
    \caption{Series solution~(\ref{newseries}) (using~(\ref{seriesfixed}) for $\tilde{\theta}(t)$) taken up to (and including) the $\bar{t}^N$ term, shown for $N$=3, 6, and 36 and compared with the numerical solution ($\bullet$) to~(\ref{pendeqnd}) for the case of  rotation ($\bar{E}>2$). The energy $\bar{E}=2.02$ is same as that in figure~\ref{fig:series}b but here the series is initialized at the top of the trajectory such that it converges. The solution is shown over 2 periods above, encompassing the solution in~figure~\ref{fig:series}b.}
    \label{fig:rotationfixed}
\end{figure}
\clearpage
\begin{figure}[h!]
    \includegraphics[width=5.1in]{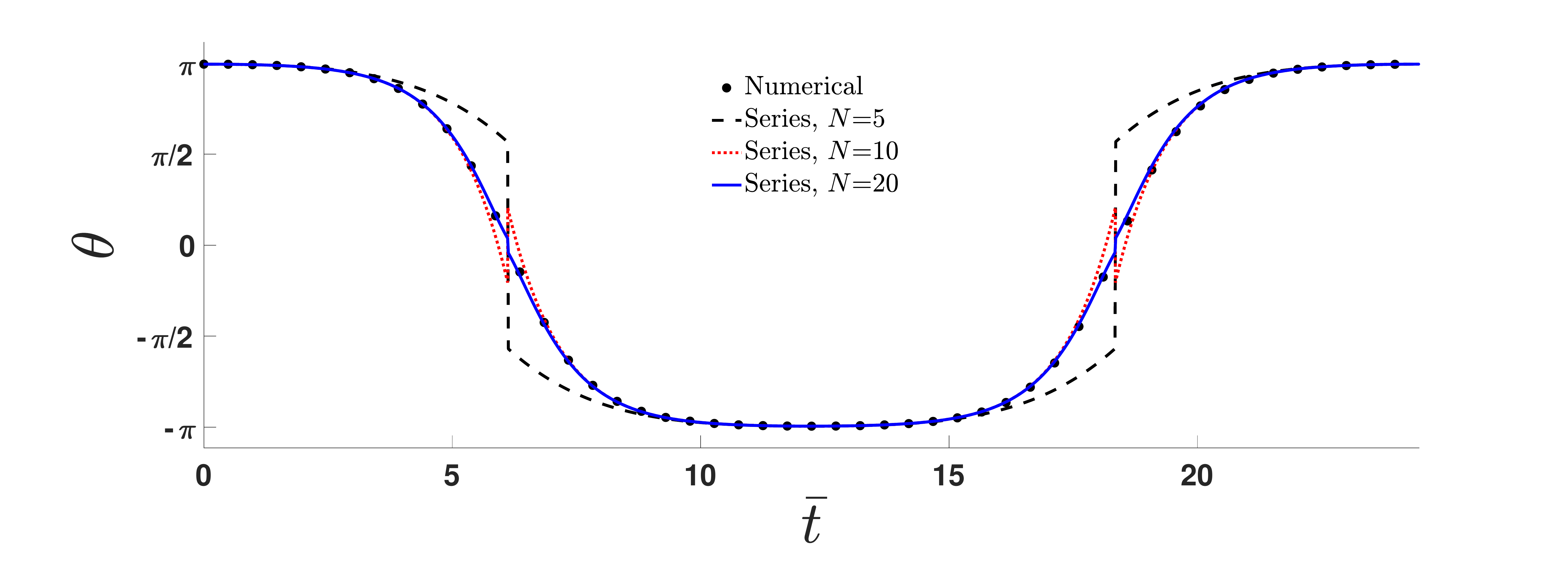}
    \caption{Series solution~(\ref{newseries}) (using~(\ref{seriesfixed}) for $\tilde{\theta}(t)$) taken up to (and including) the $\bar{t}^N$ term, shown for $N$=5, 10, and 20 and compared with the numerical solution ($\bullet$) to~(\ref{pendeqnd}) for the librating case of $\bar{E}=1.9998$. }
    \label{fig:librationseries}
\end{figure}

% \clearpage

% Although~(\ref{seriesfixed}) is an exact analytic solution to the pendulum problem with modest convergence, this convergence may be accelerated, as shown by the solid lines of figure~\ref{fig:error}. The method for convergence acceleration is given in section~\ref{sec:approximant}.

\section{Exact resummation of the series solution \label{sec:approximant}}
Although~(\ref{seriesfixed}) is an exact analytic solution to the pendulum problem, its convergence may be accelerated by a similarly exact resummation technique. The resummation is motivated by the method of asymptotic approximants~\cite{Barlow:2017} that assures that it asymptotically approaches the exact solution at each end of the domain, taken here to be $\bar{t}\in[0,T^*]$. We match the low-order behavior at $\bar{t}=T^*$ (up to the 2$^\mathrm{nd}$ derivative) and up any any desired order at $t=0$. Using the choice of initial conditions prescribed by~(\ref{seriesfixed}) to preserve analyticity, the resummation of~(\ref{seriesfixed}) is
\begin{subequations}
\begin{align} 
\nonumber
\theta(\bar{t}) = \bar{\omega}^* \left(\bar{t}-T^* \right) +  \left(\bar{t}-T^*\right)^2\sum_{n=0}^\infty\hat{a}_n\bar{t}^n,~~|\bar{t}|\le T^\star,\\ (\theta_0,~\bar{\omega}_0)=\begin{cases} \begin{rcases}
\left(\arccos(1-\bar{E}),~0\right), & \text{$\bar{E}<2$} \\
\left(\pi,~\pm\sqrt{2\bar{E}-4}\right),  & \text{$\bar{E}>2$},
 \end{rcases}, \end{cases}
\label{resummation}
\end{align} 
where
\begin{equation}
        \bar{\omega}^* \equiv \begin{cases}
   \sqrt{2\bar{E}} & \bar{E} < 2 \\
    \pm\sqrt{2\bar{E}} & \bar{E} > 2
    \end{cases}
    \label{omegab}
\end{equation}    
enforces the correct angular velocity at $\bar{t}=T^*$ and the $\pm$ symbol specifies either counterclockwise or clockwise rotation.  In~(\ref{resummation}), the coefficients $\hat{a}_n$ are related to the known coefficients $a_n$ (given by~(\ref{seriesrec})) as follows. We set~(\ref{resummation}) equal to~(\ref{sum}) to obtain
\begin{equation*}
    \sum_{n=0}^\infty a_n\bar{t}^n = \bar{\omega}^* \left(t-T^* \right) +  \left(t-T^*\right)^2\sum_{n=0}^\infty\hat{a}_nt^n, 
\end{equation*}
and then isolate the $\hat{a}_n$ series onto one side to obtain
\begin{equation*}
   \overbrace{\left[\sum_{n = 0}^\infty \left((n+1)\left(\frac{1}{T^*}\right)^{n+2}\right) \bar{t}^n\right]}^{(\bar{t}-T^*)^{-2}} \underbrace{\sum_{n=0}^\infty b_n\bar{t}^n}_{\left[-\bar{\omega}^* \left(t-T^* \right)+\displaystyle\sum_{n=0}^\infty a_n\bar{t}^n\right]} =  \sum_{n=0}^\infty\hat{a}_n\bar{t}^n.
   \label{beforeproduc}
\end{equation*}
As indicated in the above,  $(\bar{t}-T^*)^{-2}$ is expanded explicitly as an arithmetico-geometric series, and the coefficients
\begin{equation}
b_0 = a_0 + T^*\bar{\omega}^*,~b_1 = a_1 - \bar{\omega}^*,~b_{n > 1} = a_n
\label{b}
\end{equation}
are introduced to enable the left-hand side of~(\ref{beforeproduc}) to be combined into a single Cauchy product series.  The result is an expression for the $\hat{a}_n$ coefficients as 
\begin{equation}
    \hat{a}_n = \sum_{k=0}^n b_{n-k}\left((k+1)\left(\frac{1}{T^*}\right)^{k+2}\right).
    \label{approxcoeffs}
\end{equation}
\label{fullresummation}
\end{subequations}
Together,~(\ref{resummation}),~(\ref{b}), and~(\ref{approxcoeffs}) form an exact resummation of~(\ref{seriesfixed}) that explicitly enforces the asymptotic behaviors as $\bar{t}\to 0$ (from the power series~(\ref{seriesfixed})) and $\bar{t}\to T^*$, whereas~(\ref{seriesfixed}) does not have this latter feature.   We remind the reader that the restriction on initial conditions in~(\ref{fullresummation}) is not a restriction on energy; one can use~(\ref{fullresummation}) to provide convergent solutions for all librating and rotating energy states of the pendulum.  

\begin{figure}[h!]
\centering
    \includegraphics[width=3.4in]{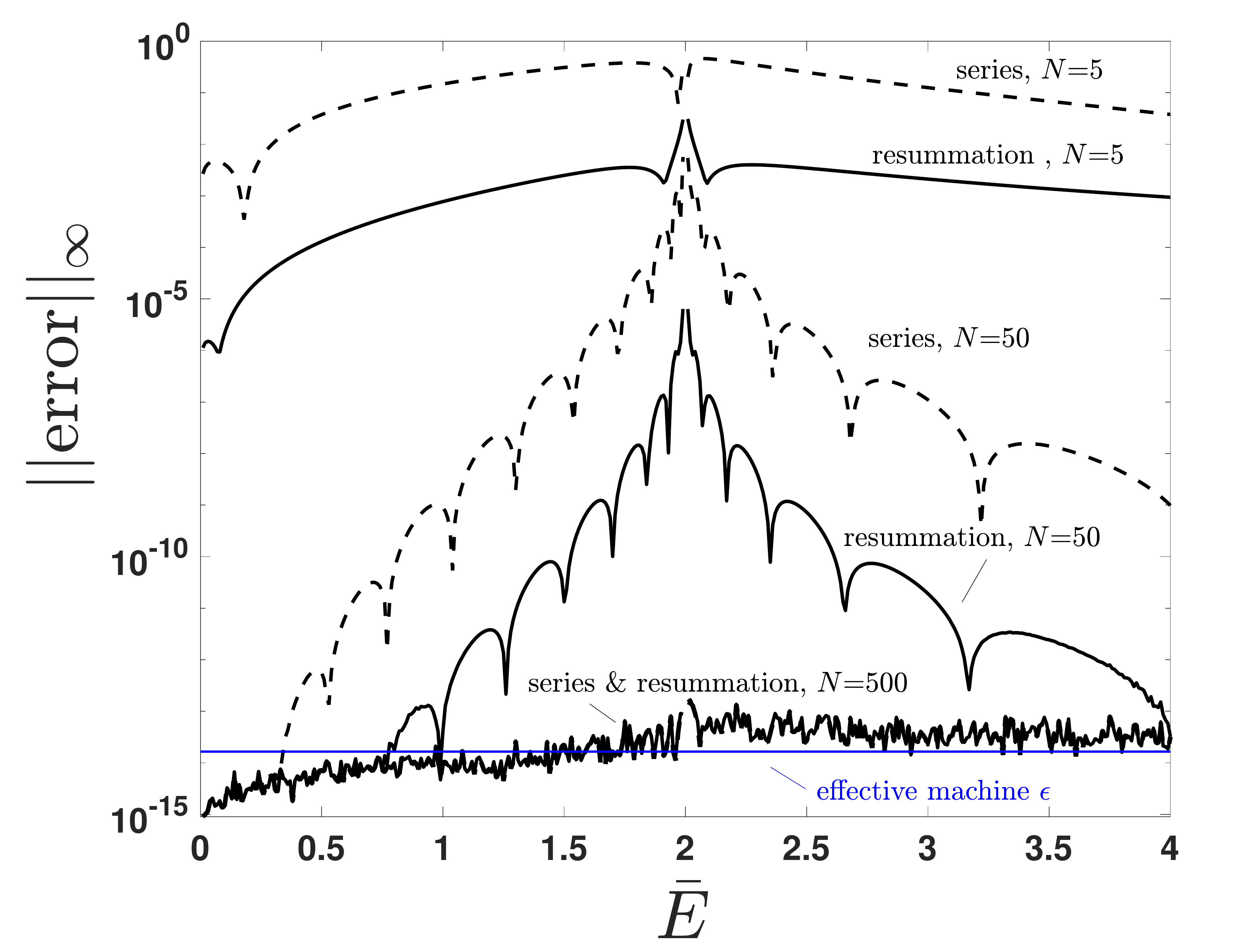}
    \caption{The error here is taken to be the difference between the numerical solution to~(\ref{pendeqnd}) and either the series solution~(\ref{seriesfixed}) or exact resummation~(\ref{fullresummation}) (as indicated in the figure) taken up to (and including) the $\bar{t}^N$ term and shown for various $N$. The infinity norm of this error is taken over the interval $\bar{t}\in[0,T^*]$ for each value of $\bar{E}$ shown, except at $\bar{E}=2$ where an exact solution exists.  The difference between the exact solution~(\ref{critical}) and the numerical solution for $\bar{E}=2$ is shown as a horizontal line, indicating the effective machine precision, as a baseline of comparison with all the other error calculations shown (for which the same numerical scheme and time-step is used).}
    \label{fig:error}
\end{figure}

Taking the resummation~(\ref{fullresummation}) as $\tilde{\theta}(t)$ in the periodically extended solution~(\ref{newseries}),  the ($N+1$)-term partial sums are shown in  figure~\ref{fig:librationapproximant}.  Inspection of figure~\ref{fig:librationseries} and figure~\ref{fig:librationapproximant},  generated for the same physical parameters, indicates that the resummation~(\ref{fullresummation}) accelerates convergence,  in that less terms are needed to achieve the same accuracy. This convergence acceleration extends to other energy values, as shown in figure~\ref{fig:error}, where the norm of the error in~(\ref{fullresummation}) is shown (by a solid line) for various ($N+1$)-term partial sums of the resummation series~(\ref{fullresummation}).
\begin{figure}[h!]
    \includegraphics[width=5.1in]{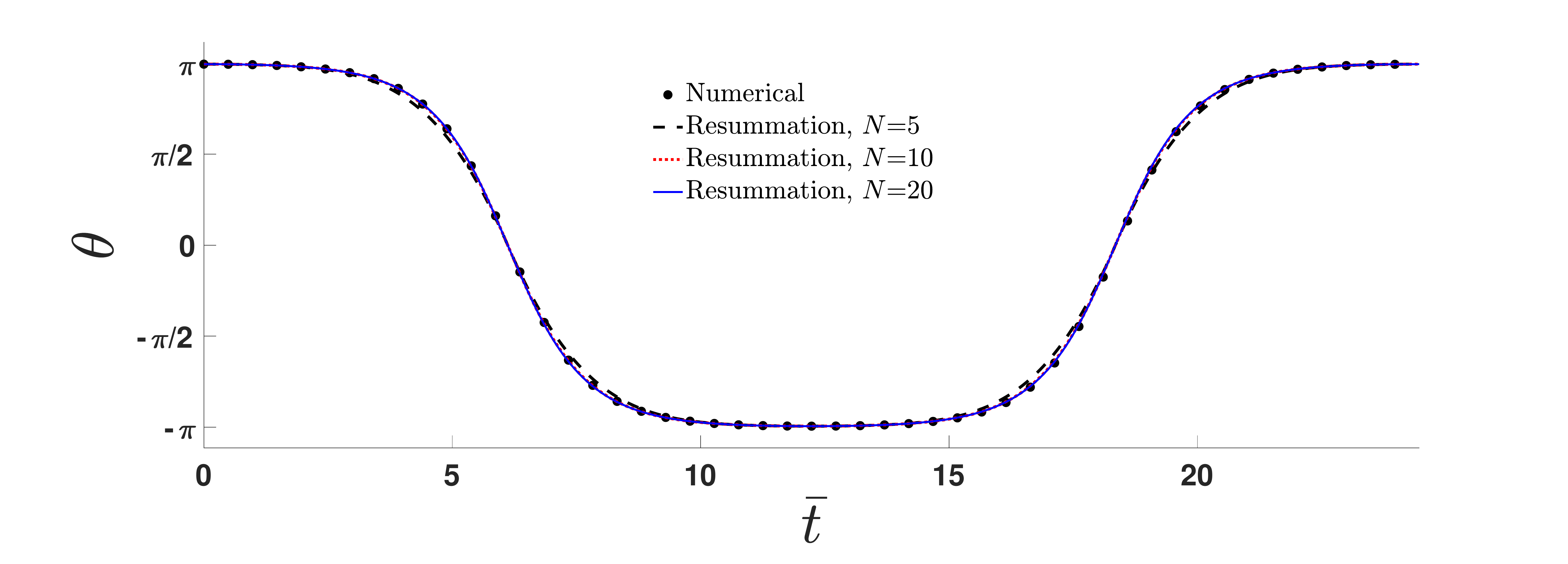}
    \caption{Exact resummation~(\ref{fullresummation}) taken up to (and including) the $\bar{t}^N$ term, shown for $N$=5, 10, and 20, periodically extended using~(\ref{newseries}), and compared with the numerical solution ($\bullet$) to~(\ref{pendeqnd}) for same conditions of figure~\ref{fig:librationseries}.  A direct comparison with figure~\ref{fig:librationseries} shows the convergence advantage imparted by the resummation series~(\ref{fullresummation}) over the original series~(\ref{seriesfixed}).}
    \label{fig:librationapproximant}
\end{figure}
Note that further manipulations, shown in Appendix~\ref{sec:A1},  may be used to increase computational speed for any finite truncation of the resummation series~(\ref{fullresummation}).

\section{Exact resummation of the pendulum period
\label{sec:period}}
In both the series solution~(\ref{seriesfixed}) and its resummation~(\ref{fullresummation}), the series converge for $\bar{t}\in[0,T^*]$, and $T^*$ explicitly appears in both~(\ref{seriesfixed}) and the periodic extension formulae given in~(\ref{newseries}); thus, the value of $T^*$ (and by extension, the value of $T$) must be determined accurately. As can be seen in~(\ref{T}) and~(\ref{Tstar}), in every case of period motion, $T^*$ is written in terms of the complete elliptic integral of the first kind, $K(k)$, where $k$ depends on the energy of the pendulum.  Thus we are faced with computing the period either by numerical integration, using an approximation (of which there are many~\cite{Hinrichsen}), or using a truncation of the long-known infinite series representation of $K(k)$~\cite{Byrd,Abramowitz}, namely:
\begin{equation}
     K(k) = \frac{\pi}{2}\sum_{n=0}^\infty \left(\frac{(2n)!}{n!^2}\right)^2 16^{-n}   k^{2n},~~k^2<1. 
     \label{Kold}
\end{equation}
In~(\ref{Kold}), note that the $k^2<1$ constraint is always satisfied for the pendulum problem (see~(\ref{T}), for example).   The error in calculating $T^*$ using~(\ref{Kold}) in~(\ref{Tstar}) is shown by the dashed lines in figure~\ref{fig:perioderror}, plotted versus $\bar{E}$ for partial sums up to and including the $k^{2N}$ term (i.e, $N+1$ terms in~(\ref{Kold})).  Note that the error increases as $\bar{E}\to 2$ where the period becomes infinite, as it corresponds to the separatrix in figure~\ref{fig:surf}. 

\begin{figure}
\centering
    \includegraphics[width=3.5in]{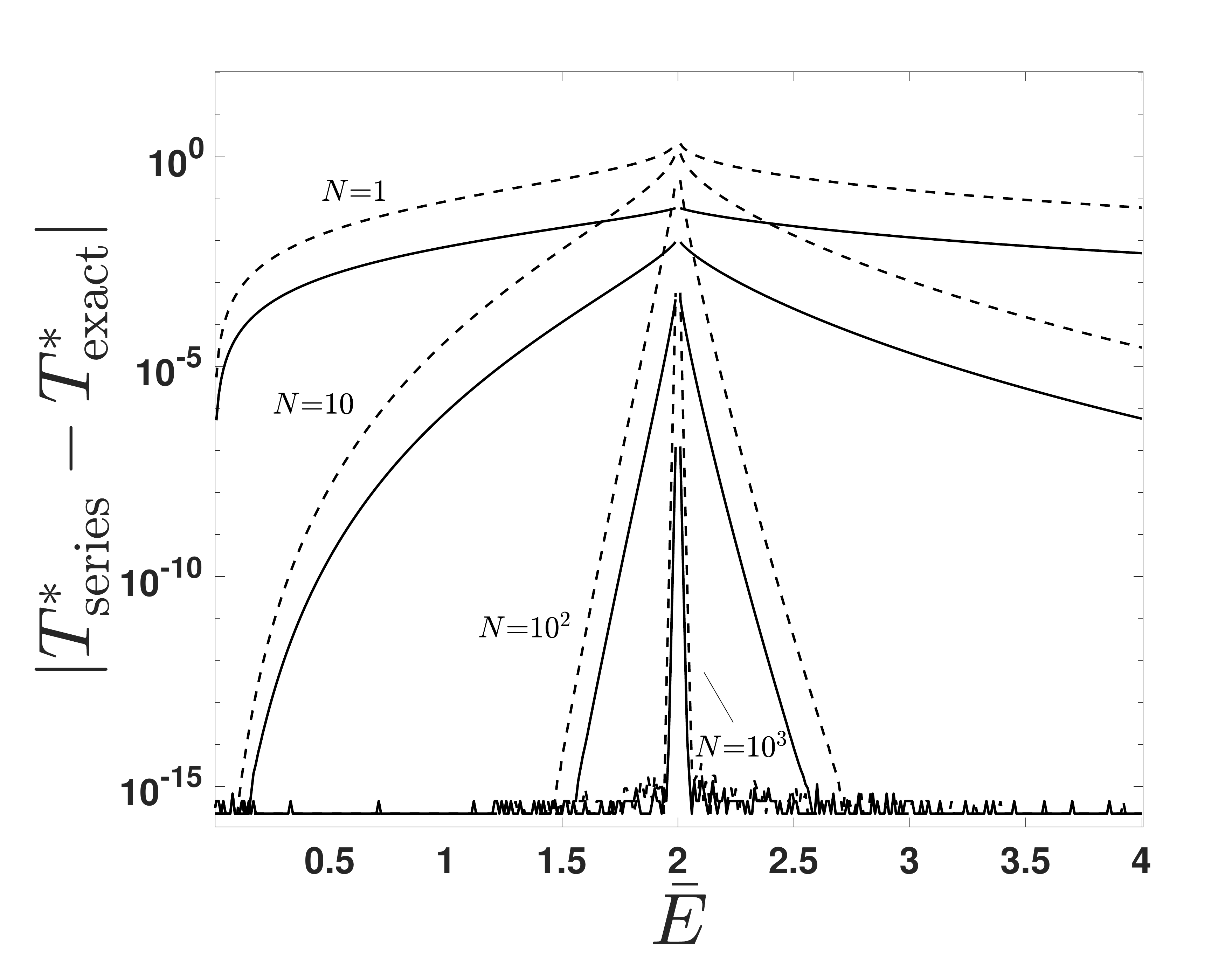}
    \caption{The error is shown for evaluating a pendulum's effective period $T^*$ (and by extension, actual period) given by~(\ref{Tstar}).  In evaluating $T^*_\mathrm{series}$, the elliptic integral $K(k)$ given by~(\ref{K}), is taken to be either the original series~(\ref{Kold}) (dashed line) or its resummation~(\ref{Knew}) (solid line) taken up to (and including) the $k^{2N}$ term.  The exact solution, denoted here as $T^*_\mathrm{exact}$, uses a numerical evaluation of~(\ref{K}) that has converged to within machine precision.}
    \label{fig:perioderror}
\end{figure}

Figure~\ref{fig:period} shows results close to $\bar{E}=2$ (above and below) that show the effect on the periodic extension of solution~(\ref{fullresummation}) via~(\ref{newseries}) when~(\ref{Kold}) is used to evaluate $T^*$ (see dashed curves).  In the figure, a truncation of $N$=100 is used in~(\ref{Kold}) in order in illustrate that one would need to take more than 100 terms in the series~(\ref{Kold}) before a solution using this value of $T^*$ can be considered accurate. 

\begin{figure}
\centering
    \includegraphics[width=3.5in]{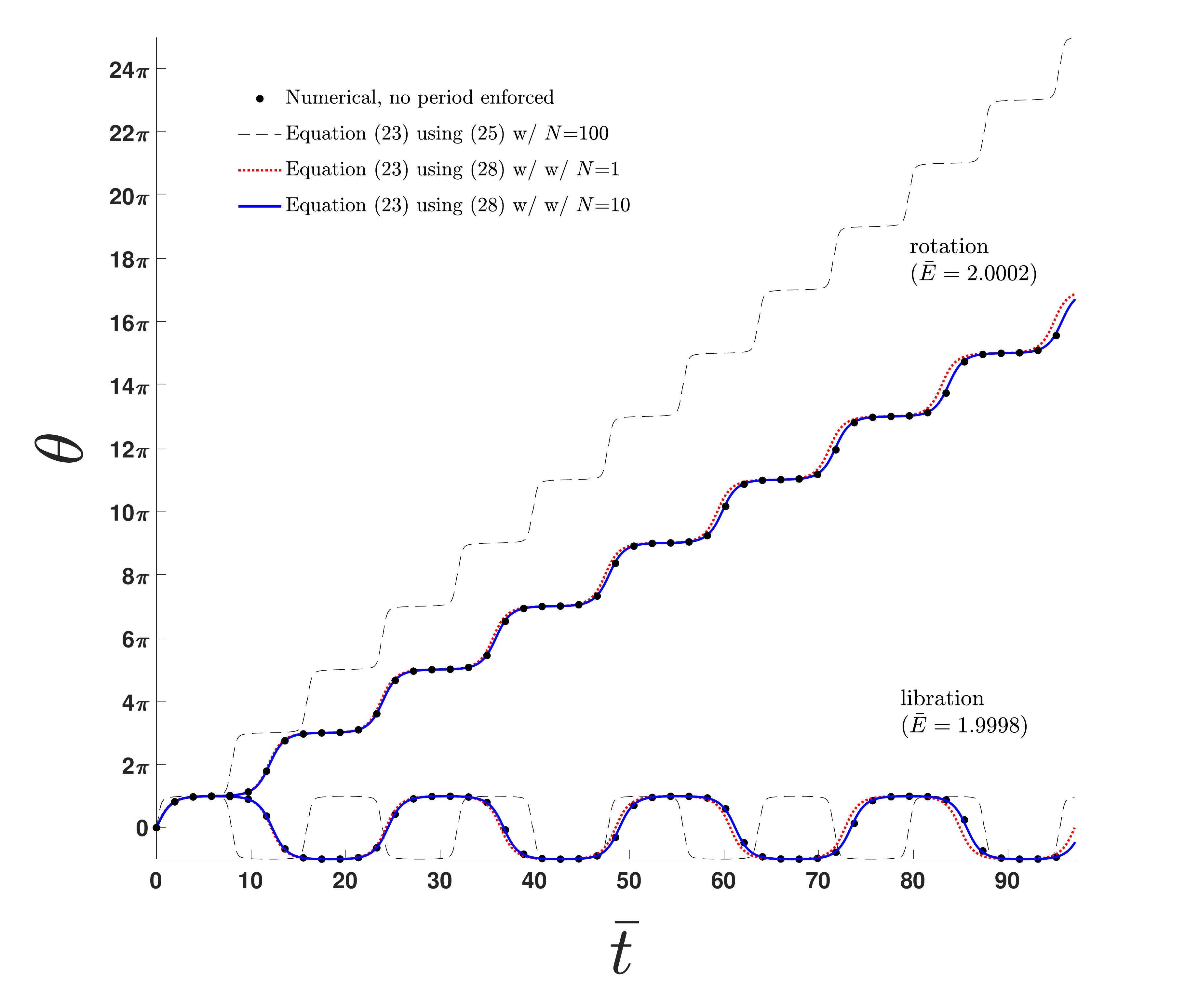}
    \caption{ Here, we demonstrate the effect of period accuracy on the periodically extending~(\ref{fullresummation}) via~(\ref{newseries}) (taken up to 20$^\mathrm{th}$ order), shown for librating and rotating cases near $\bar{E}=2$ as indicated in the figure. The numerical solution is given by $\bullet$s. The long-dashed curves show a solution that uses a period computed using~(\ref{T}) with~(\ref{Kold}) and letting $N=100$.  The short-dashed and solid curves show a solution that uses a period computed using~(\ref{T}) with~(\ref{Knew}) and letting $N=1$ and $N=10$, respectively.}
    \label{fig:period}
\end{figure}

% \TODO{}\comment{used $\theta_0=0$ and $\bar{\omega}_0= 1.9999$ and $2.0001$}
% Thus we note that the method of asymptotic approximants can also produce approximations to be used in the quarter period that are competitive with the state-of-the-art amongst analytic expressions for the period that can be readily and efficiently computed to arbitrarily high order. 
In order to make practical use of~(\ref{Kold}), we accelerate its convergence using a resummation approach.  To do so, we  begin with an equivalent definition to~(\ref{K}) for the complete elliptic integral of the first kind~\cite{Byrd}: 
\begin{equation}
   K(k) = \int_0^1 \left[(1-v^2)(1-k^2v^2)\right]^{-1/2}\ \text{d}v. 
   \label{K2}
\end{equation}
We follow a similar approach to that used in~\cite{Xue2012} to strip off the singular behavior of~(\ref{K2}) by first rewriting the integral as $$ K(k) = \underbrace{\int_0^1 \frac{\sqrt{1-k^2v^2}-\sqrt{1-v^2}}{(1-k^2v^2)\sqrt{1-v^2}}\ \text{d}v}_{I_1}  + \underbrace{\int_0^1 \frac{1}{1-k^2v^2}\ \text{d}v}_{I_2}   $$ so that $I_2$  naturally handles the apparent integrand singularity at $k=1$. Using partial fractions, this integral may be evaluated exactly as $$I_2 =  \frac{1}{2k} \ln\left(\frac{1+k}{1-k}\right), $$ and Taylor expanded as  \begin{equation}
    I_2 =\frac{1}{2k} \ln\left(\frac{1+k}{1-k}\right) =  \frac{\operatorname{arctanh}(k)}{k}= \sum_{n=0}^\infty \frac{1}{2n+1}k^{2n}.
    \label{zero}
\end{equation}   Thus, we may use~(\ref{zero}) to add 0 to~(\ref{Kold}) in such a way that the singular behavior (responsible for slow convergence) is exactly summed, leading to an improved series for $K(k)$, given as \begin{equation}
   K(k) = \sum_{n=0}^\infty \left. \vphantom{\frac{\pi}{2}\left(\frac{(2n)!}{n!^2}\right)^2 16^{-n}} \right[ \frac{\pi}{2}\left(\frac{(2n)!}{n!^2}\right)^2 16^{-n} - \frac{1}{2n+1}\left. \vphantom{\frac{\pi}{2}\left(\frac{(2n)!}{n!^2}\right)^2 16^{-n}} \right]  k^{2n} + \frac{1}{2k} \ln\left(\frac{1+k}{1-k}\right),~~k^2<1.
   \label{Knew}
\end{equation}   
  When used to compute $T^*$ (or $T$) in~(\ref{Tstar}) (or~(\ref{T})) the convergence of the resummation~(\ref{Knew}) is significantly improved, compared with that of~(\ref{Kold}), as seen in figure~\ref{fig:perioderror} for a range of $\bar{E}$ and in figure~\ref{fig:period} for $\bar{E}$ close to 2.

  %Note that there are several other infinite expansions for $K(k)$~\cite{Abramowitz,Byrd}; we use~(\ref{Knew}) here for its simplicity and accuracy (as seen in figures~\ref{fig:perioderror} and~\ref{fig:period}).    

% Empirically, the convergence of this approximant is good, and significantly better than expressions with the same asymptotic behavior that come less naturally from the integral representation. The central binomial coefficient in equation \ref{} can be written in terms of double factorials in ways that can be computed recursively for more efficient dynamic programming representations of the above. 

\section{Conclusions \label{sec:conclusions}}
The exact series solution to the pendulum problem, whether librating or rotating, attains optimal convergence when the initial conditions are chosen to be at the top of the pendulum's trajectory, allowing the series to converge beyond the bottom of the trajectory. Through symmetry, this convergent solution is then used to construct the solution for all time.  Taking inspiration from the method of asymptotic approximants, the exact power series solution is re-summed to accelerate convergence, based on the asymptotic behavior at the bottom of the trajectory.  The resulting exact power series representation is shown to be more accurate than the original series for any given series truncation and pendulum energy.

We close this paper with some general comments about the methodology used here as well as the solution itself.  When finding the series solution to a nonlinear ODE, there are (at least) 3 scenarios: (i) the radius of convergence is larger than the physical domain, (ii) the radius of convergence is precisely the physical domain length (e.g., series for Rayleigh's equation~\cite{Harkin}), and (iii) the radius of convergence lies within the physical domain (e.g, Blasius series~\cite{Blasius}).   The exact pendulum solution provided in this paper falls into case (i), once the initial conditions are chosen with care.  Here our analysis benefited from the properties of the Jacobi elliptic functions whose singularity locations are known.  In case (ii), a singularity is known at the boundary of the domain and this is the closest singularity in the complex plane to the expansion point.  Case (iii) arises from spontaneous singularities and is commonplace in the solution of nonlinear ODEs.  Note that even when case (iii) occurs, resummation methods such as asymptotic approximants may be used to analytically continue a series over the whole domain to obtain a useful convergent solution (see~\cite{Barlow:2017,FS2020,SIR2020,SEIR2020}).

In contrast to the problem examined here, there is often no pre-knowledge of singularity locations in the solution of nonlinear ODEs. For such cases, one can find the radius of convergence for a series around a given expansion point numerically using a Domb-Sykes plot~\cite{DombSykes} or root test.
For autonomous equations such as this one, then, an optimal expansion point would need to be deduced numerically; whether or not this expansion point leads an exact solution depends on the final radius of convergence achieved.  Note, however, that even in cases where the radius of convergence is not larger than the domain of interest, approximant techniques may be utilized to analytically continue such series and make them useful.  Coupled with approximants (when needed), the power series solution method is a prevailing technique that can lead to efficient, useful, and many times exact solutions to nonlinear ODEs.

\appendix

\section{Efficient computation of~(\ref{fullresummation}) \label{sec:A1}}
Recognizing that the exact resummation of the pendulum series, given by~(\ref{fullresummation}), will be implemented in any code as a partial sum, we provide an efficient means of calculating the deflection angle. One may observe that the truncation of~(\ref{resummation}), denoted here as $\theta_N(t)$, has the following structure
\begin{subequations}
\begin{equation}
    \theta_{N}(t) = \underbrace{\bar{\omega}^*\left(t-T^*\right)}_{\text{degree 
1 polynomial}} + \underbrace{\left(t-T^*\right)^2}_{\text{degree 2 polynomial}}\underbrace{\hphantom{\int}\sum_{n=0}^N \hat{a}_nt^n}_{\text{degee $N$ polynomial}}.
\label{truncation}
\end{equation}
Since the approximant given by~(\ref{truncation}) is just a degree $N+2$ polynomial, we may rewrite it as
\begin{equation}
    \theta_{N}(t) = \sum_{n=0}^N a_n t^n + \alpha t^{N+1} + \beta t^{N+2}.
    \label{truncation2}
\end{equation}
Upon equating the above to the original series~(\ref{sum}), one finds that the $a_n$ coefficients are precisely the original series coefficients given in~(\ref{pendseries}). Hence, in order to obtain the approximant in the form of~(\ref{truncation2}), we first determine the $a_n$ from the formulae give in~(\ref{pendseries}), and then use those to compute $\alpha$ and $\beta$ so that the first and second derivatives of~(\ref{truncation}) and~(\ref{truncation2}) match at $\bar{t}=T^*$. This leads to a system of two linear equations in $\alpha$ and $\beta$, whose solution is

\begin{align}
    \alpha &= -(N+2)\sigma_N({T^*}) {T^*}^{-N-1} - \left(\bar{\omega}^* - \sigma_N'({T^*})\right) {T^*}^{-N} \\ 
    \beta &= (N+1) \sigma_N({T^*}) {T^*}^{-N-2} + \left(\bar{\omega}^* - \sigma_N'({T^*})\right){T^*}^{-N-1}
\end{align}

where 

\begin{equation}
    \sigma_{N}(t)\equiv \sum_{n=0}^N a_n t^{N}.
\end{equation}
\label{efficient}
\end{subequations}

For comparison, the entire process to compute the resummation coefficients (for a given truncation $N$) via the formulae given in~(\ref{fullresummation}) (including the computation of series coefficients $a_n$ from~(\ref{pendseries})) requires two quadratic time steps. The first step is the recursion in~(\ref{pendseries}) used to get the series solution, and the second is Cauchy's product rule, used in~(\ref{fullresummation}). The method given in this appendix \-- recasting~(\ref{fullresummation}) as~(\ref{efficient}) \-- removes the latter of those quadratic time steps, which leads to a substantial reduction in computational time \-- such as when the series is re-used for different values of $\bar{E}$ to construct the solution surface in figure~\ref{fig:surf}.

\bibliographystyle{unsrt}
\bibliography{pendulum}

\end{document}